\documentclass[11pt,reqno]{amsart} 
\usepackage{graphicx}
\usepackage{amssymb,amsthm,listings,enumerate,geometry,setspace,float,csquotes,mathtools,url}
\usepackage{etoolbox}
\usepackage[labelformat=parens,labelfont={}]{subcaption}
\usepackage{hyperref}
\usepackage{pdfpages,xspace,xparse}
\usepackage{import,pst-plot}
\usepackage[numbers,sort]{natbib}
\usepackage{bm}

\DeclareRobustCommand{\bms}[1]{\bm{#1}}
\newcommand{\MATLAB}{\textsc{Matlab}\xspace}

\theoremstyle{definition}
\newtheorem{remark}{Remark}[section]
\numberwithin{equation}{section}
\numberwithin{figure}{section}

\makeatletter
\def\@settitle{\begin{center}%
  \baselineskip14\p@\relax
  \@title
  \end{center}%
}

\title[Num. solns. of lin. SD problems on half-line using UTM]{\Large{The numerical solutions of linear semi-discrete evolution problems on the half-line using the Unified Transform Method}\vspace{-15pt}}

\author[J. Cisneros \& B. Deconinck]{}

\begin{document}
\maketitle
\begin{center}
	Jorge Cisneros$^{1}$ \& Bernard Deconinck$^{2}$ \\[5pt]
	Department of Applied Mathematics \\
	University of Washington\\
	Seattle, WA 98195-2420\\ 
	$^{1}$\texttt{jorgec5@uw.edu}, $^{2}$\texttt{deconinc@uw.edu}\\[5pt]
	\today
\end{center}

\begin{abstract}
We discuss a semi-discrete analogue of the Unified Transform Method, introduced by A. S. Fokas, to solve initial-boundary-value problems for linear evolution partial differential equations of constant coefficients. The semi-discrete method is applied to various spacial discretizations of several first and second-order linear equations on the half-line $x \geq 0$, producing the exact solution for the semi-discrete problem, given appropriate initial and boundary data. We additionally show how the Unified Transform Method treats derivative boundary conditions and ghost points introduced by the choice of discretization stencil. We consider the continuum limit of the semi-discrete solutions and provide several numerical examples.\\[7pt]
\fontsize{8}{9.6}\selectfont{Keywords: Unified Transform Method, semi-discrete linear problem, half-line, ghost points, continuum limit}
\end{abstract}

\section{Introduction}

	Consider the numerical solution of the $N^{\text{th}}$-order quasilinear partial differential equation (PDE)
	\begin{equation}
		q_t = c \, q_{Nx} + F\left(q,q_x,\ldots,q_{(N-1)x} \right),\quad c \in \mathbb{C} \, \text{\textbackslash} \, \{ 0 \},
		\label{ibvp_eq}
	\end{equation} 
	on the half-line, $x \in (0,\infty)$ or on the finite interval, $x \in (0,L)$ with $L > 0$. The solution of \eqref{ibvp_eq} is uniquely determined if we additionally prescribe an initial condition $q(x,0) = \phi(x)$ and the correct number of boundary conditions that are compatible at $x = t = 0$ with sufficient smoothness and decay. For the half-line problem, we additionally impose initial and boundary conditions so that $q(x,t) \rightarrow 0$ uniformly as $x \rightarrow \infty$.
	
\begin{figure}[b]
	\begin{center}
  		\includegraphics[width=0.5\linewidth]{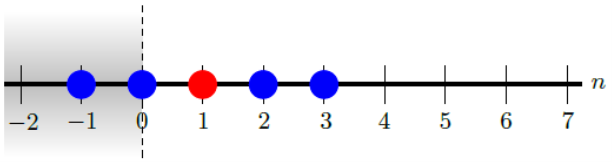}
		\caption{A stencil that requires information at a ghost point.}
		\label{ghost_point}
	\end{center}
\end{figure}

	Arguably the most intuitive approach in solving finite-interval initial-boundary value problems (IBVPs) is through the implementation of a finite-difference scheme on a discrete grid with points $x_n \equiv n \Delta x$ and $t_j \equiv j \Delta t$. Directly applying such schemes, especially those with high-order spatial stencils, introduces the dependence on grid points outside of the domain, known as \textit{ghost points}, see Figure \ref{ghost_point}. This embeds a discrepancy into the numerical methodology, since these points originate from the choice of spatial stencil and not from the original IBVP itself. Note that with periodic boundary conditions, the issue of ghost points never arises. 

	The choice of information at ghost points can destabilize numerical methods that are shown to be stable in the full-line or periodic problem via von Neumann analysis \cite{randy,dale}. These heuristic methodologies do not easily transfer when treating PDEs with higher-order derivatives, and the general rules for examining stability in the presence of boundary data are not well developed \cite{gustafsson,cheema_thesis,coco,weno,black_hole}. Incorporating boundary conditions correctly and addressing ghost points is a non-trivial numerical issue \cite{iserles,strikwerda,trefethen}. 

	Our approach to tackle this problem is set up by the ideas of operator splitting and the implementation of split-step methods. Higher-order derivatives tend to require higher-order stencils, so we consider the class of semi-discretized PDEs where the most nonlocal stencil is applied to the linear term $c \, q_{Nx}$. Hence, the lower-order problem $q_t = F\left(q,q_x,\ldots,q_{(N-1)x} \right)$ of \eqref{ibvp_eq} can be approached using established split-step techniques, while the linear problem $q_t = c \, q_{Nx}$ and the ghost points that arise require special attention. In this paper, we treat $q_t = c \, q_{Nx}$ using the Unified Transform Method on the semi-discrete $(n,t)$-plane, with an eye towards split step methods for \eqref{ibvp_eq} to be explored in a future paper. The same treatment for finite-interval problems will also appear in a subsequent paper \cite{SDUTM_FI}.


\section{The Continuous Unified Transform Method}\label{cont_UTM_steps}

	The Unified Transform Method (UTM) or Method of Fokas provides a powerful approach to solve evolution IBVPs, including all those with linear, constant-coefficient PDEs and some integrable nonlinear PDEs. The UTM was introduced by A. S. Fokas in 1997 for the purpose of generalizing the method of inverse scattering to IBVPs on the half-line and on a finite interval \cite{fokas_paper,fokas_collab,fokas_book}. 

	The UTM generates an explicit analytical solution for $q(x,t)$, with the solution written in terms of integrals along paths in the complex plane of a spectral parameter $k \in \mathbb{C}$. Through parametrization of the contours, these explicit solutions can be numerically evaluated \cite{flyer,xin_paper}. The application of the UTM is systematic, regardless of the types of boundary data, \textit{e.g.}, nonhomogeneous Dirichlet, Neumann, Robin, etc. conditions. This is one reason the UTM is more general and effective than standard methods for evolution IBVPs. Further, the method demonstrates how many and which types of boundary conditions result in a well-posed IBVP, depending on the order $N$ of the PDE \cite{fokas_book}.
	
	For either half-line or finite-interval IBVPs, the UTM is applied algorithmically using the following steps \cite{bernard_fokas}:
\begin{enumerate}
	\item Rewrite the PDE in divergence form, depending on a spectral parameter $k$, to obtain the local relation and the dispersion relation $\tilde{W}(k)$,\label{utm_1}
	\item Integrate over the $(x,T)$-domain to obtain the global relation, \label{utm_2}
	\item Invert the global relation to obtain a representation of the solution depending on known and unknown boundary data, \label{utm_3}
	\item Determine symmetries $\tilde{\nu}_j(k)$ of $\tilde{W}(k)$, \label{utm_4}
	\item Determine where in $\mathbb{C}$ the global relations evaluated at $\tilde{\nu}_j(k)$ are valid, \label{utm_5}
	\item If necessary, deform integral paths involving boundary terms appropriately, \label{utm_6}
	\item Solve for unknown boundary data using the global relations evaluated at $\tilde{\nu}_j(k)$, and \label{utm_7}
	\item Check that integral terms involving $\hat{q}(\tilde{\nu}_j,T)$ vanish, resulting in a solution representation. \label{utm_8}
\end{enumerate}
	Although the calculations within each step are more intricate for higher-order problems and their boundary conditions, the UTM ultimately solves an IBVP by solving a set of algebraic equations involving the dispersion relation and its symmetries.
	
	For the finite-difference evaluations of nonlinear IBVPs, we can apply the UTM to the semi-discrete problem. A method-of-lines formulation allows the UTM to address ghost points directly by providing an analytical solution to the linear semi-discrete IBVP.


\section{Semi-Discrete UTM: Notation and Definitions}

	The UTM has received a lot of attention for continuous IBVPs, but not nearly as much for semi-discrete ones, \textit{i.e.}, discretized in space $x_n = n h$, but continuous in time $t$. Biondini \& Hwang \cite{gino_nls}, Biondini \& Wang \cite{gino_main}, and Moon \& Hwang \cite{sd_finite} study semi-discrete problems in the context of the semi-discrete UTM (SD-UTM), but from the perspective of a purely semi-discrete problem on closed contours with discretized Lax pairs and no variable mesh spacing $h$. While \cite{gino_nls} focuses primarily on the linear and nonlinear Schr\"{o}dinger difference equations, the theory for semi-discrete problems is presented via examples in \cite{gino_main} for half-line IBVPs and in \cite{sd_finite} for finite-interval IBVPs. Minimal discussion on the continuum limit for the SD-UTM is presented in \cite{sd_finite}. 
	
	Our goal is to further develop the SD-UTM to help solve IBVPs by addressing complications that arise with ghost points. Within the split-step method, we want to apply the SD-UTM to the linear problem, whether that includes solving the linear semi-discrete IBVP entirely or only applying the method at the ghost points generated from the choice of spatial stencil in the finite-difference scheme. As we will see next, the SD-UTM formulas for the semi-discrete $q_n(T)$ are simpler than those from the continuous UTM, but further approximations are needed in order to efficiently implement them into a split-step method, see Section \ref{small_time_sec}.  
	
	In what follows, we present the SD-UTM through examples of several linear semi-discretized IBVPs on the half-line, with an explicit mesh parameter $h \ll 1$. For each section, the first few concrete examples are followed by higher-order discretizations where ghost points arise. We follow a similar procedure to Steps \eqref{utm_1}-\eqref{utm_8}, outlining any changes as we move through the examples. We use the shift operator $\Delta Q_n = Q_{n+1} - Q_n$, which effectively replaces the spatial derivative with a forward difference. For IBVPs with Dirichlet boundary conditions, the Fourier transform pair can be written as
	\begin{subequations}
		\begin{equation}
			\hat{q}(k,t) = h \sum_{n = 1}^{\infty} e^{-ik nh} q_n(t) , \quad\quad \text{Im}(k) \leq 0,
			\label{fourier_cont_SD}
		\end{equation}
		\begin{equation}
			q_n(t) = \frac{1}{2 \pi} \int_{-\pi/h}^{\pi/h} e^{iknh} \hat{q}(k,t) \, dk, \quad\quad k \in \mathbb{C}.
			\label{inv_fourier_cont_SD}
		\end{equation}
		\label{fourier_transforms_SD}
	\end{subequations}
	If a Dirichlet boundary is not given, then \eqref{fourier_cont_SD} starts at $n = 0$. For half-line IBVPs, we require $q_n \in l^1(\mathbb{N})$, the space of absolutely summable sequences, ensuring that $\hat{q}(k, t)$ is bounded for all $k \in \mathbb{C}$ with $\text{Im}(k) \leq 0$. Additionally, we define the time transform of spatial \textit{nodes} at and near the $n = 0$ boundary, including ghost points: 
	\begin{equation}
		f_j(W,T) = \int_{0}^T e^{Wt} q_j(t) \, dt, \quad\quad k \in \mathbb{C},
		\label{F}
	\end{equation}
	with a finite $T > 0$ and semi-discrete dispersion relation $W(k)$. 
	
	We briefly discuss the difficulties of numerically computing the solution to half-line IBVPs via finite-difference schemes. Conventionally, we truncate the half-line problem $x \in [0,\infty)$ to a finite-interval problem $x \in [0, M]$, where $M \in \mathbb{R}$ is a large positive constant, so that the artificial numerical boundary is far from the domain of interest $x \in [0,L]$. At the artificial boundary $x = M$, we can apply, say, decaying boundary conditions that are compatible with the given initial condition. Now, the half-line IBVP is recast as a finite-interval problem and the usual finite-difference tactics can be applied. This approach heavily relies on $M \gg L$, so that contributions from the artificial boundary do not interfere with the window of interest. For dispersive problems, the effect of a tail slowly approaches zero, and $M$ might have to be prohibitively large, increasing the computational cost to produce an accurate solution. Below, we do not compare the semi-discrete UTM solutions with traditional windowing finite-difference methods. 
	
	Embedded within the following examples, we compute the SD-UTM solutions within the window of interest $x \in (0,1]$ if a Dirichlet boundary condition is specified or $x \in [0,1]$ if a Neumann boundary condition is given. The solutions are implemented in \MATLAB using built-in functions, such as the vectorized \texttt{integral()}. To reduce computation time, we analytically evaluate sums, like those defining the forward discrete Fourier transform, and integrals when possible. In addition, all IBVPs have initial and boundary conditions matching at $(x,t) = (0,0)$. For the second-order problems, the exact solutions to the continuous problems are given in terms of error functions, which are well optimized for numerical evaluations.


\section{Advection Equations}\label{advec_eqns_sec}

To start, we discuss advection equations in some detail, as a way to demonstrate the UTM applied to semi-discrete problems. At the same time, this will allow us to fix notation and to illustrate the types of numerical experiments we use throughout the paper.


\subsection{Forward Discretization of $\bms{q_t = c\, q_{x}}$}\label{advec_forward_halfline}
\begin{figure}[tb]
		\begin{center}
			\includegraphics[width=0.4\linewidth]{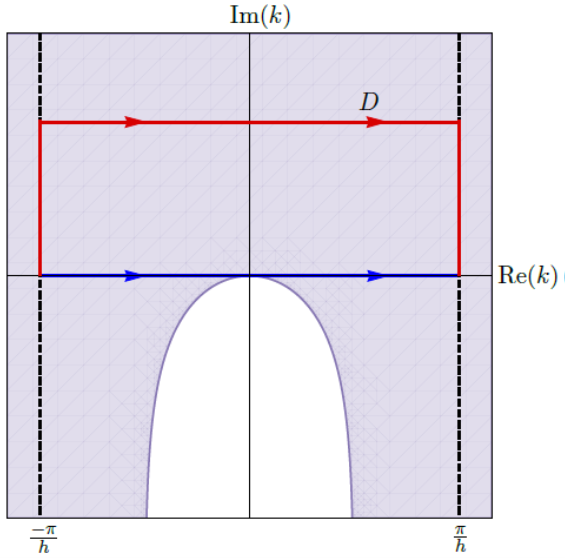}
			\caption{The shaded regions depict where $\text{Re}(-W) \leq 0$ and $e^{-WT}$ is bounded with the dispersion relation \eqref{W_advec1_forward}.}
			\label{advec_forward_W}
		\end{center}
	\end{figure}

	We start with the continuous problem on the half-line for the advection equation $q_t = c\, q_{x}$ with wave-speed $c > 0$:
		\begin{equation}\begin{dcases}
		q_t = c \, q_{x},& x > 0,\, t > 0, \\
		q(x,0) = \phi(x),& x > 0.
		\label{advec1_prob}
	\end{dcases}\end{equation}
	For well posedness, the IBVP requires only the initial condition and no boundary data. Since information travels from right to left, it is well known that the forward discretization of $q_{x}(x,t)$ together with a forward discretization of $q_t(x,t)$ is a ``natural'' discretization, known as the upwind method. Such a method performs well for this advection equation with periodic boundary conditions or on the whole real line with $\lim_{x \rightarrow \pm \infty} q(x,t) = 0$. Let us implement this forward spacial discretization. We consider
	\begin{align}
		\dot{q}_n(t) = c\,\frac{q_{n+1}(t) - q_{n}(t)}{h},
		\label{advec1_forward}
	\end{align}
	followed by using the semi-discrete version of the UTM to \textit{exactly} solve this system of ODEs, instead of a time-stepping method to \textit{approximately} (because of the time discretization) solve the system. As in the continuous UTM, the local relation is determined by writing the problem into its divergence form. For this semi-discrete problem, we replace $\partial_x$ with the shift operator $\Delta Q_n$, and \eqref{advec1_forward} is rewritten as
	\begin{equation}
		\partial_t \left(e^{-iknh} e^{Wt} q_n \right) = \frac{c}{h}\Delta \left(e^{-ik(n-1)h} e^{Wt} q_{n} \right),
		\label{LR_advec1_forward}
	\end{equation}
	with dispersion relation 
	\begin{equation}
		W(k) = c\, \frac{1 - e^{ikh}}{h}.
		\label{W_advec1_forward}
	\end{equation}
	The symmetries of a dispersion relation are those transformations $k \rightarrow \nu(k)$ that leave $W(k)$ invariant, \textit{i.e.}, $W(\nu) = W(k)$. Here, \eqref{W_advec1_forward} only has the trivial symmetry $\nu_0(k) = k$ up to periodic copies due to the complex exponential. From the local relation \eqref{LR_advec1_forward}, we obtain the global relation by taking a time transform over $t \in [0,T]$ and an infinite sum from $n = 0$ (because $q_0(t)$ is not known):
	\begin{align}
		&&\hspace{-55pt} \sum_{n=0}^{\infty} h \int_0^T \left[ \partial_t \left(e^{-iknh} e^{Wt} q_n \right) - \frac{c}{h}\Delta \left(e^{-ik(n-1)h} q_{n} \right)e^{Wt} \right] dt &= 0 \notag \\
		\Rightarrow \,&&\hspace{-55pt} \sum_{n=0}^{\infty} h \left[ e^{-iknh} e^{WT} q_n(T) - e^{-iknh} q_n(0) - \frac{c}{h}\Delta \left(e^{-ik(n-1)h}  f_{n} \right) \right] &= 0 \notag \\
		\Rightarrow \,&&\hspace{-55pt} e^{WT} \hat{q}(k,T) - \hat{q}(k,0) + c e^{ikh} f_{0} &= 0,
		\label{GR_advec1_forward}
	\end{align}
	valid for $\text{Im}(k) \leq 0$ due to the discrete Fourier transform terms. Solving for $\hat{q}(k,T)$ and inverting the inverse transform,
	\begin{align}\begin{split}
		q_n(T) &= \frac{1}{2\pi} \int_{-\pi/h}^{\pi/h} e^{iknh} e^{-WT} \hat{q}(k,0)\,dk - \frac{c}{2\pi} \int_{-\pi/h}^{\pi/h} e^{ik(n+1)h} e^{-WT} f_0\,dk.
		\label{soln1_advec_forward}
	\end{split}\end{align}
	 The integrand in the first term is defined for $\text{Im}(k) \leq 0$, while the integrand in the second term is defined for all $k \in \mathbb{C}$. We refer to the expression above as the ``solution,'' since $f_0(W,T)$ in the second integral term is unknown. For all $n \in \mathbb{N}$, $e^{ik(n+1)h}$ decays in the upper half-plane and $e^{-WT}$ is bounded in the shaded regions, including the  boundary, of Figure \ref{advec_forward_W}. The shaded region denotes where $\text{Re}(-W) \leq 0$. Figure \ref{advec_forward_W} also shows the integration path for ``solution'' \eqref{soln1_advec_forward} from $-\pi/h$ to $\pi/h$ on the real line.
	Note that the sign of $c$ is essential in determining the location of the region of exponential growth of the integrand, \textit{i.e.}, the white region in Figure \ref{advec_forward_W}.
	
	We use two approaches to address the unwanted boundary integral term in ``solution'' \eqref{soln1_advec_forward}. The first approach is more straightforward, but is not as general as the second approach. In both, we substitute the definition of $f_0(W,T)$ in order to collect the $k$ dependence:
	\begin{align*}
		\frac{c}{2\pi} \int_{-\pi/h}^{\pi/h} e^{ik(n+1)h} e^{-WT} f_0\,dk &= \frac{c}{2\pi} \int_{-\pi/h}^{\pi/h} e^{ik(n+1)h} e^{-WT} \left[\int_{0}^T e^{Wt} q_0(t)\,dt \right]\,dk = \int_{0}^T A(n,T-t) q_0(t)\,dt,
	\end{align*} 
	with $T - t > 0$ and
	$$A(n,T) = \frac{c}{2\pi} \int_{-\pi/h}^{\pi/h} e^{ik(n+1)h} e^{-WT} \,dk.$$
	
	\begin{enumerate}[(i)]
	\item{
	The first approach uses the transformation $z = e^{ikh}$:
	\begin{align*}
		A(n,T) = \frac{c}{2 \pi i h} \oint_{|z|=1} z^{n} \exp\left[-\left(\frac{1-z}{h}\right)(T-t)\right] \,dz = 0,
	\end{align*}
	by analyticity of the integrand for all $n$. Hence,
	$$\frac{c}{2\pi} \int_{-\pi/h}^{\pi/h} e^{ik(n+1)h} e^{-WT} f_0\,dk = 0.$$
	}
	\item{
	The second approach deforms the integration path of $A(n,T)$ away from the real line. Consider $R > 0$. We define the line segment
	$$D = \left\{k \in \mathbb{C} \, \Big| \, \frac{-\pi}{h} \leq \text{Re}(k) \leq \frac{\pi}{h} \, \text{ and } \, \text{Im}(k) = R \right\}$$
	with left-to-right orientation. Thus, $D$ is a horizontal straight-line path above the real line, from $k = \tfrac{-\pi}{h} + i R$ to $k = \tfrac{\pi}{h} + i R$. Next, we introduce a closed contour that consists of four straight segments: the original real-line path, the new path $D$, and two vertical segments that connect the endpoints of the real-line path with those of $D$, as illustrated in Figure \ref{advec_forward_W}. The contribution to the integral from these vertical paths cancel due to periodicity. Hence,
	$$A(n,T) = \frac{c}{2\pi} \int_{D} e^{ik(n+1)h} e^{-W(T-t)} \,dk,$$
	by Cauchy's Theorem. Taking $R \rightarrow \infty$ implies taking $\text{Im}(k) \rightarrow \infty$ in the integrand. Because of the exponential decay above the real line, $A(n,T) = 0$ and
	$$\frac{c}{2\pi} \int_{-\pi/h}^{\pi/h} e^{ik(n+1)h} e^{-WT} f_0\,dk = 0.$$
	}
	\end{enumerate}
	It follows that the solution to the half-line IBVP with the forward discretization \eqref{advec1_forward} depends only on the initial condition:
	\begin{align}
		q_n(T) &= \frac{1}{2\pi} \int_{-\pi/h}^{\pi/h} e^{iknh} e^{-WT} \hat{q}(k,0)\,dk.
		\label{soln_advec_forward}
	\end{align}
	
	For reference, we solve the IBVP \eqref{advec1_prob} using the continuous UTM, following the Steps \eqref{utm_1} -- \eqref{utm_8} from Section \ref{cont_UTM_steps}. Briefly, we find the dispersion relation $\tilde{W}(k) = - c i k$, with only the trivial symmetry $\nu_0(k) = k$, and the global relation 
	$$\hat{q}(k,0) - e^{\tilde{W}T} \hat{q}(k,T) + c F_0 = 0,\quad \text{Im}(k) \leq 0,$$
	where
	$$\hat{q}(k,t) = \int_{0}^{\infty} e^{-ikx} q(x,t) \, dx ,\quad \text{Im}(k) \leq 0,$$
	and
	$$F_j(\tilde{W},T) = \int_{0}^T e^{\tilde{W}t} \left. \frac{\partial^j q}{d x^j} \right|_{x = 0} \, dt, \quad k \in \mathbb{C}.$$
	After inverting the transform and showing there is no dependence on $F_0(\tilde{W},T)$, the solution representation is 
	\begin{align}
		q(x,T) &= \frac{1}{2 \pi} \int_{-\infty}^{\infty} e^{ikx} e^{-\tilde{W}T} \hat{q}(k,0) \, dk.
		\label{soln_advec1_cont}
	\end{align}
	Taking the limit as $h \rightarrow 0$ of \eqref{soln_advec_forward}, we recover \eqref{soln_advec1_cont} from the continuous problem, where the limits of integration approach $\pm \infty$ at rate $1/h$. Also, $\lim_{h \rightarrow 0}W(k) = - c i k = \tilde{W}$.

\begin{figure}[tb]	
		\raggedleft
		\begin{subfigure}[t]{.45\textwidth}
			\centering
  			\includegraphics[width=0.95\linewidth]{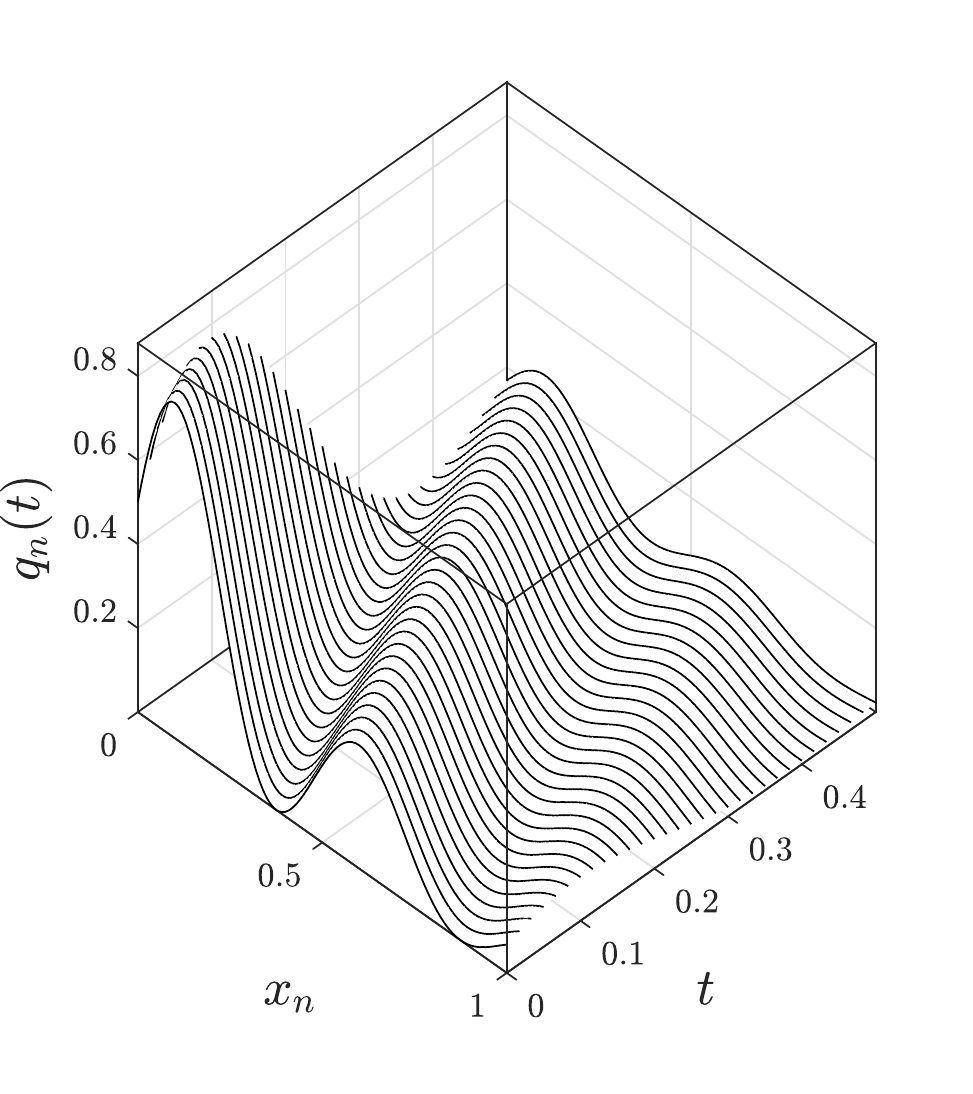}
  			\caption{}
  			\label{advec1_UTM1_HL_solnplot}
		\end{subfigure}\hfill 
		\begin{subfigure}[t]{.45\textwidth}
			\centering
  			\includegraphics[width=1\linewidth]{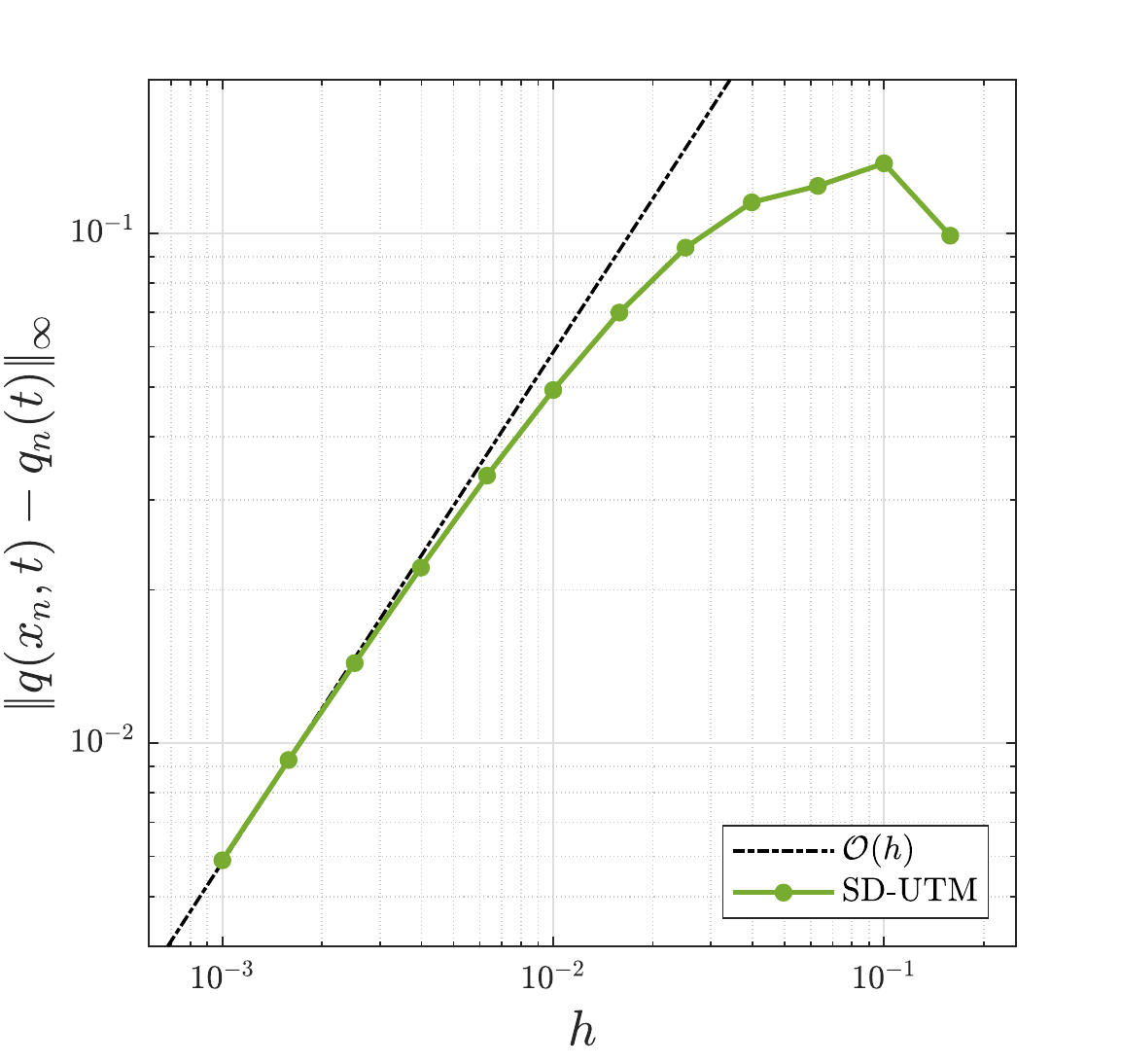}
  			\caption{}
  			\label{}
		\end{subfigure}	
		\caption{(a) The semi-discrete solution \eqref{soln_advec_forward} evaluated at various $t$ with $h = 0.01$. (b) Error plot of the semi-discrete solution \eqref{soln_advec_forward} relative to the exact solution as $h \rightarrow 0$ with $t = 0.5$.}
		\label{advec1_UTM1_HL}
	\end{figure}

	As an explicit example, we compute the numerical solution to the IBVP 
	\begin{equation}
	\begin{dcases}
		q_t = q_{x},& x > 0,\, t > 0, \\
		q(x,0) = \phi(x) = \tfrac{1}{2}\left[ e^{-2x} \left( \sin(4 \pi x) + 1\right)\right],& x > 0.
	\end{dcases}
	\label{advec1_numerical1_HL}
	\end{equation}
	The exact (continuous) solution is given by $q(x,T) = \phi(x+T)$, while the semi-discrete solution is obtained from the representation \eqref{soln_advec_forward} with the standard forward discretization stencil. Figure \ref{advec1_UTM1_HL} shows the semi-discrete solution $q_n(t)$ (left panel) and a log-log error plot (right panel) of the $\infty$-norm of $q_n(0.5)-q(x_n, 0.5)$, as a function of $h$.
	
	From the stencil \eqref{advec1_forward}, we know \eqref{soln_advec_forward} is a first-order accurate approximation to the solution $q(x,T)$ of the IBVP \eqref{advec1_prob}. We can reveal more information about the behavior and structure of this approximate solution by determining its modified equation \cite{randy}. Suppose $q_n(T)$ exactly solves a PDE with dependent variable $p(x,T)$, such that $q_n(T) \equiv p\left(x_n,T\right)$. Substituting this assumption into the forward stencil \eqref{advec1_forward} and Taylor-series expanding terms gives
	\begin{align*}
		\dot{q}_n(t) &= \frac{c}{h} \left[ q_{n+1}(t) - q_{n}(t) \right] \\
		\Rightarrow \quad p_t(x_n,t) &= \frac{c}{h} \left[ p(x_n+h,t) - p(x_n,t) \right] \\
		&= \frac{c}{h} \left[ p(x_n,t) + p_x(x_n,t) h + \frac{p_{xx}(x_n,t)}{2!} h^2 + \frac{p_{xxx}(x_n,t)}{3!} h^3 + \mathcal{O}\left(h^4\right) - p(x_n,t) \right] \\
		p_t &= c\, p_x + \frac{c\,p_{xx}}{2} h + \frac{c\, p_{xxx}}{6} h^2 + \mathcal{O}\left(h^3\right).
	\end{align*}
	Keeping up to the $\mathcal{O}(h)$ term, we find that \eqref{soln_advec_forward} is a second-order accurate solution approximation to the advection-diffusion PDE
	\begin{equation}
		p_t = c \,p_x + \frac{c\,h}{2} p_{xx},
		\label{advec_forward_modified_eqn}
	\end{equation}
	so we expect solution profiles of \eqref{soln_advec_forward} to travel at the correct speed $c$, while dissipating in time. Since $c >0$, the diffusion coefficient $c\,h/2$ is positive. If we allow $c < 0$ or if we apply the same forward stencil to the PDE $q_t = - a \,q_x$ with $a >0$, we obtain a similar convection-diffusion modified PDE like above, except with a negative diffusion coefficient that presents an ill-posed problem with exponentially growing solutions. The solution plot \ref{advec1_UTM1_HL_solnplot} displays the expected shift to the left as time progresses. For this advection equation, the solution approaches zero as $t \rightarrow \infty$, because $q(x,0)$ decays as $x \rightarrow \infty$. The dissipative behavior from the modified PDE \eqref{advec_forward_modified_eqn} does not appear to be troublesome. The error plot displays $\mathcal{O}(h)$ convergence as $h \rightarrow 0$. 
	
	\begin{remark}
	All forward discretizations produce $f_j(W,T)$ terms with a coefficient $C_j\,e^{i \gamma_j kh}$ for some $C_j \in \mathbb{C}$ and $\gamma_j \in \mathbb{N}$. Coupled with polynomial dispersion relations $W(z)$, we can remove all integral terms containing any $f_j(W,T)$ from ``solutions'' using the steps above. Thus, if we solve the IBVP \eqref{advec1_prob} without boundary conditions using a purely forward higher-order stencil, we find \eqref{soln_advec_forward} as the solution, except with a different dispersion relation $W(k)$.
	\end{remark}


\subsection{Backward Discretization of $\bms{q_t = -c\,q_{x}}$}\label{advec2_backward_halfline}
\begin{figure}[tb]
		\begin{center}
			\includegraphics[width=0.4\linewidth]{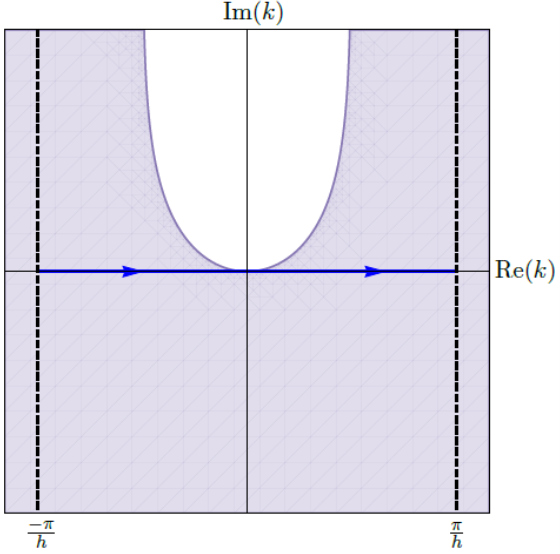}
			\caption{The shaded regions depict where $\text{Re}(-W) \leq 0$ and $e^{-WT}$ is bounded with the dispersion relation \eqref{W_advec2_backward}.}
			\label{advec_backward_W}
		\end{center}
	\end{figure}

	Next, we consider
		\begin{equation}\begin{dcases}
		q_t = -c\,q_{x},& x > 0,\, t > 0, \\
		q(x,0) = \phi(x),& x > 0,\\
		q(0,t) = u(t),& t > 0,
		\label{advec2_prob}
	\end{dcases}\end{equation}
	with $c > 0$. For well posedness, the IBVP requires a Dirichlet boundary condition at $x = 0$. Since a forward discretization \eqref{advec1_forward} for $q_t = c\,q_x$ was appropriate, we now apply a backward discretization to the spacial derivative $q_x$, resulting in
	\begin{align}
		\dot{q}_n(t) = - c\left(\frac{q_{n}(t) - q_{n-1}(t)}{h}\right) = c\,\frac{q_{n-1}(t) - q_{n}(t)}{h}.
		\label{advec2_backward}
	\end{align}
	Following similar steps as before, the local relation is 
	\begin{align}
		\partial_t \left(e^{-iknh} e^{Wt} q_n \right) &= \frac{-c}{h}\Delta \left(e^{-iknh} e^{Wt} q_{n-1} \right),
		\label{LR_advec2_backward}
	\end{align}
	with dispersion relation 
	\begin{equation}
		W(k) = c\,\frac{1 - e^{-ikh}}{h}.
		\label{W_advec2_backward}
	\end{equation}
	As before, we only have the trivial symmetry $\nu_0(k) = k$, up to periodic copies. This time, the IBVP \eqref{advec2_prob} contains a Dirichlet boundary condition, providing information at $n = 0$, so we define the forward transform as
	$$\hat{q}(k,t) =  h \sum_{n=1}^\infty e^{-iknh} q_n(t),$$
	starting at $n = 1$. To obtain the global relation, we proceed as before:
	\begin{align}
		\sum_{n=1}^{\infty} h \int_0^T \left[ \partial_t \left(e^{-iknh} e^{Wt} q_n \right) + \frac{c}{h}\Delta \left(e^{-iknh} e^{Wt} q_{n-1} \right) \right] dt &= 0 \notag \\
		\Rightarrow \hspace{25pt} e^{WT} \hat{q}(k,T) - \hat{q}(k,0) - c e^{-ikh} f_{0} &= 0,
		\label{GR_advec2_backward_transform}
	\end{align}
	valid for $\text{Im}(k) \leq 0$. Solving for $\hat{q}(k,T)$ and inverting, we obtain
	\begin{align}\begin{split}
		q_n(T) &= \frac{1}{2\pi} \int_{-\pi/h}^{\pi/h} e^{iknh} e^{-WT} \hat{q}(k,0)\,dk + \frac{c}{2\pi} \int_{-\pi/h}^{\pi/h} e^{ik(n-1)h} e^{-WT} f_{0}\,dk.
		\label{soln_advec2_backward}
	\end{split}\end{align}
	Since $e^{-WT}$ grows in the upper-half plane, see Figure \ref{advec_backward_W}, we cannot remove the dependence on $f_0(W,T)$ and, hence, \eqref{soln_advec2_backward} is the actual solution to the backward-discretized IBVP \eqref{advec2_prob} with a given Dirichlet boundary condition.
	
	Similar to the IBVP \eqref{advec1_prob}, we solve \eqref{advec2_prob} using the continuous UTM. We find the dispersion relation $\tilde{W}(k) = c i k$, with the trivial symmetry $\nu_0(k) = k$, and  
		\begin{align}
		\hspace{-25pt}q(x,T) &= \frac{1}{2 \pi} \int_{-\infty}^{\infty} e^{ikx} e^{-\tilde{W}T} \hat{q}(k,0) \, dk + \frac{c}{2 \pi} \int_{-\infty}^{\infty} e^{ikx} e^{-\tilde{W}T} F_0 \, dk.
		\label{soln_advec2_cont}
	\end{align}
	It is clear that \eqref{soln_advec2_backward} converges to the continuous solution \eqref{soln_advec2_cont}, where
	$$\lim_{h \rightarrow 0} f_j (W,T) =\lim_{h \rightarrow 0} \int_{0}^T e^{Wt} q(jh,t) \, dt = \int_{0}^T e^{\tilde{W}t} q(0,t) \, dt = F_0(W,T),$$
	for any fixed $j$, $\lim_{h \rightarrow 0}W(k) = \tilde{W}(k)$, and $\lim_{h\rightarrow0} e^{ik(n-1)h} = e^{ikx}$ with $\lim_{h \rightarrow 0} n h = x$.
	
\begin{figure}[tb]	
		\raggedleft
		\begin{subfigure}[t]{.45\textwidth}
			\centering
  			\includegraphics[width=0.95\linewidth]{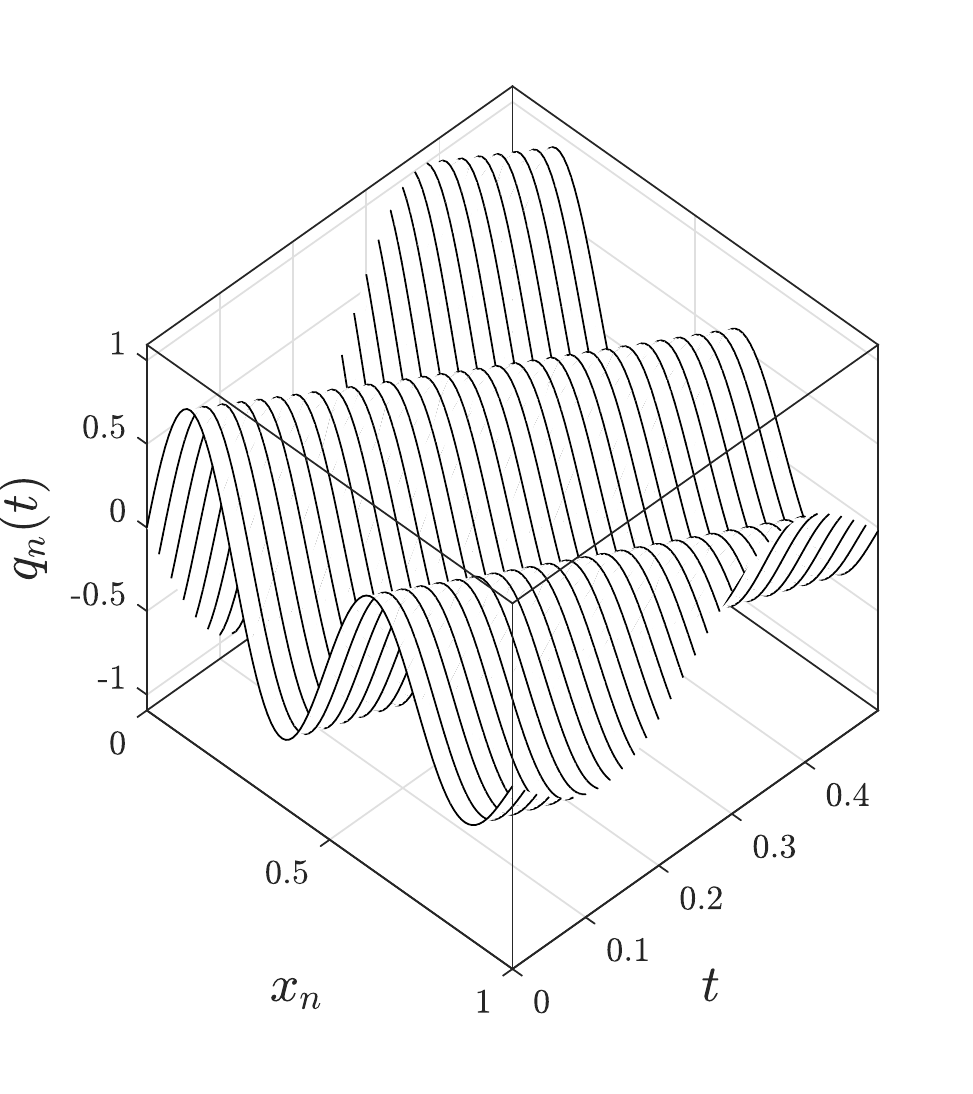}
  			\caption{}
  			\label{}
		\end{subfigure}\hfill 
		\begin{subfigure}[t]{.45\textwidth}
			\centering
  			\includegraphics[width=1\linewidth]{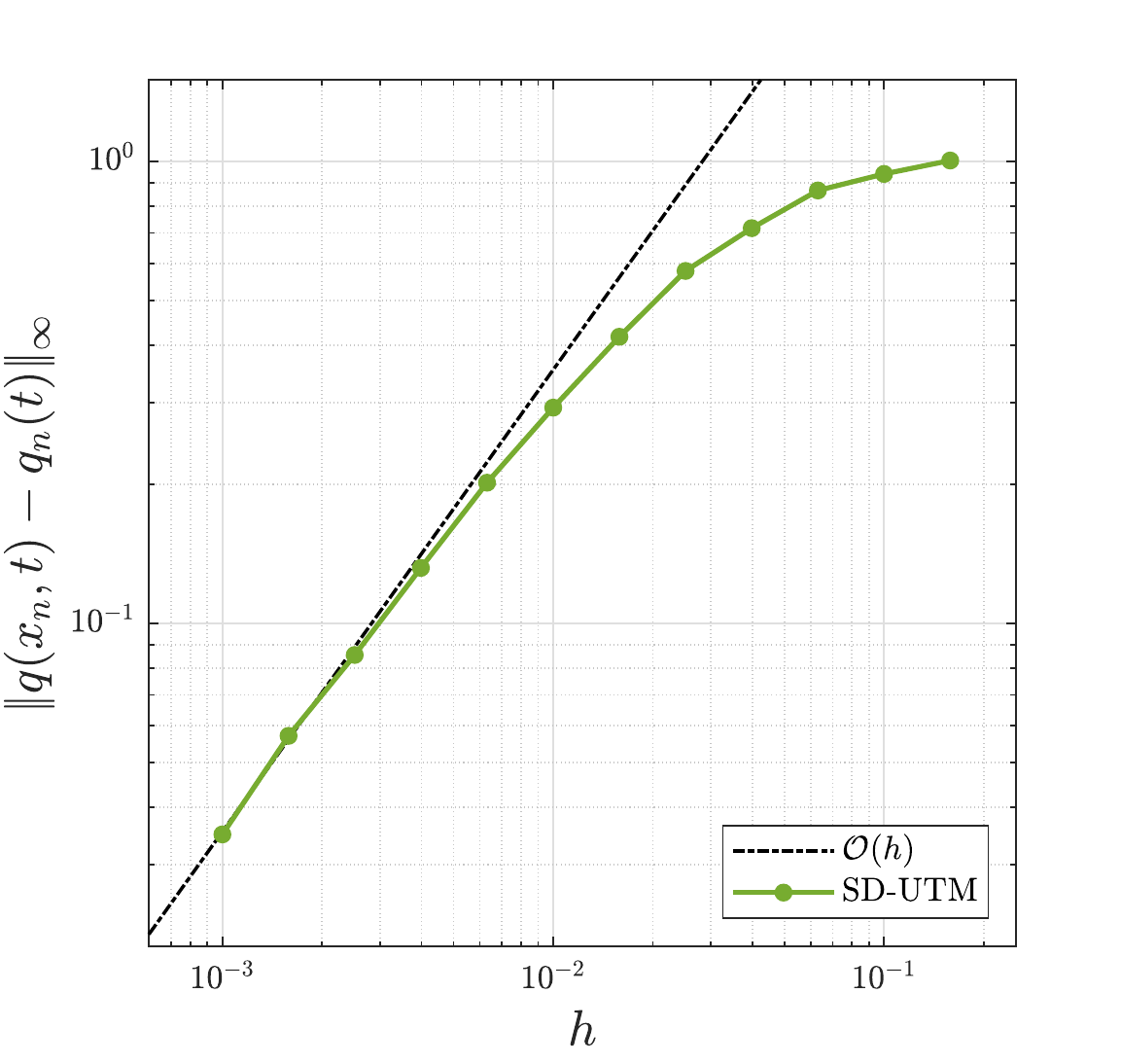}
  			\caption{}
  			\label{}
		\end{subfigure}	
		\caption{(a) The semi-discrete solution \eqref{soln_advec2_backward} evaluated at various $t$ with $h = 0.01$. (b) Error plot of the semi-discrete solution \eqref{soln_advec2_backward} relative to the exact solution as $h \rightarrow 0$ with $t = 0.5$.}
		\label{advec2_UTM2_HL}
	\end{figure}

\begin{figure}[tb]
	\raggedleft
	\begin{subfigure}{.33333\textwidth}
		\centering
  		\hspace*{-55pt}\includegraphics[width=1.1\linewidth]{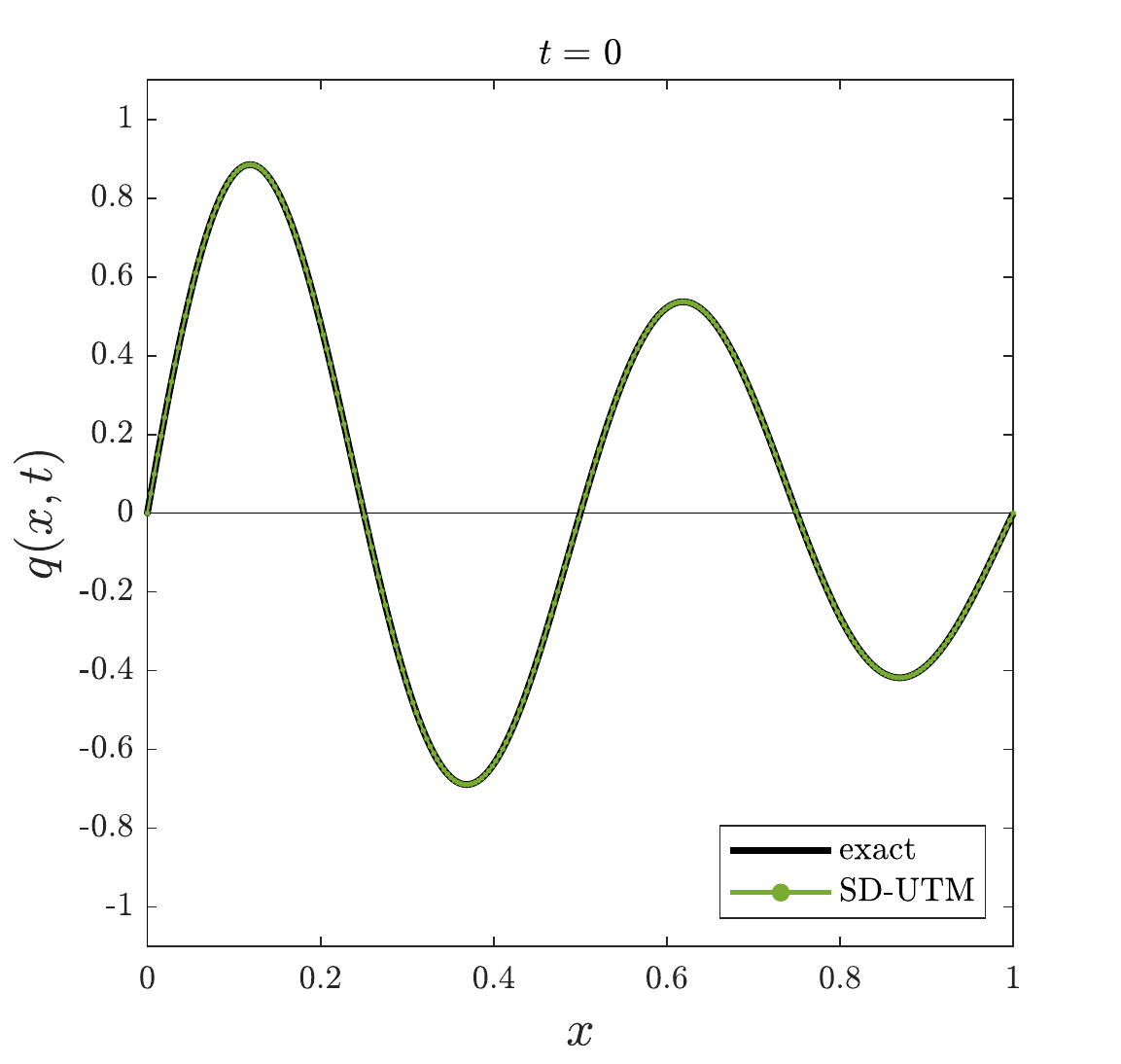}
	\end{subfigure}%
	\begin{subfigure}{.33333\textwidth}
		\centering
  		\hspace*{-5pt}\includegraphics[width=1.1\linewidth]{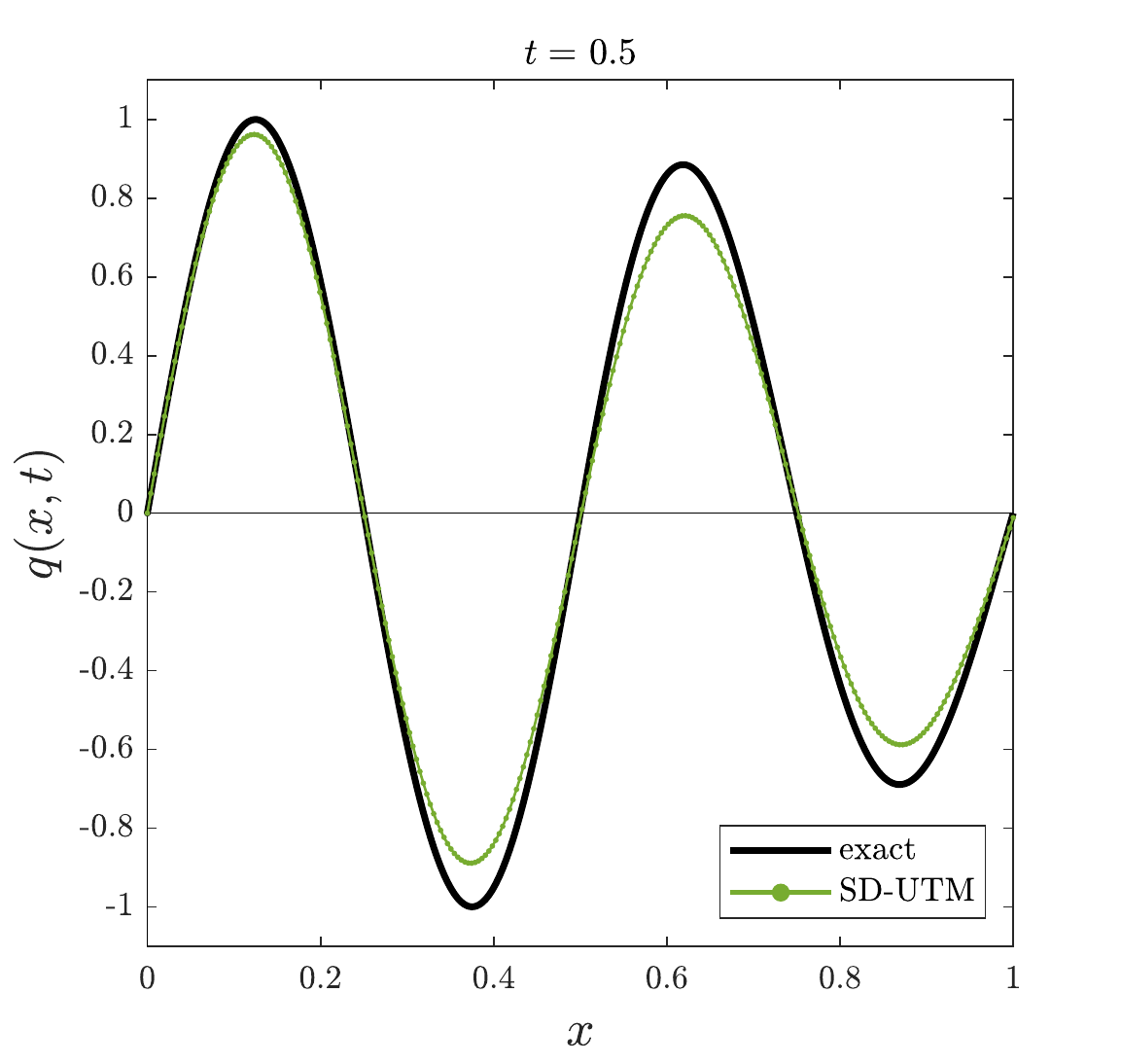}
	\end{subfigure}%
	\begin{subfigure}{.33333\textwidth}
		\centering
  		\hspace*{25pt}\includegraphics[width=1.1\linewidth]{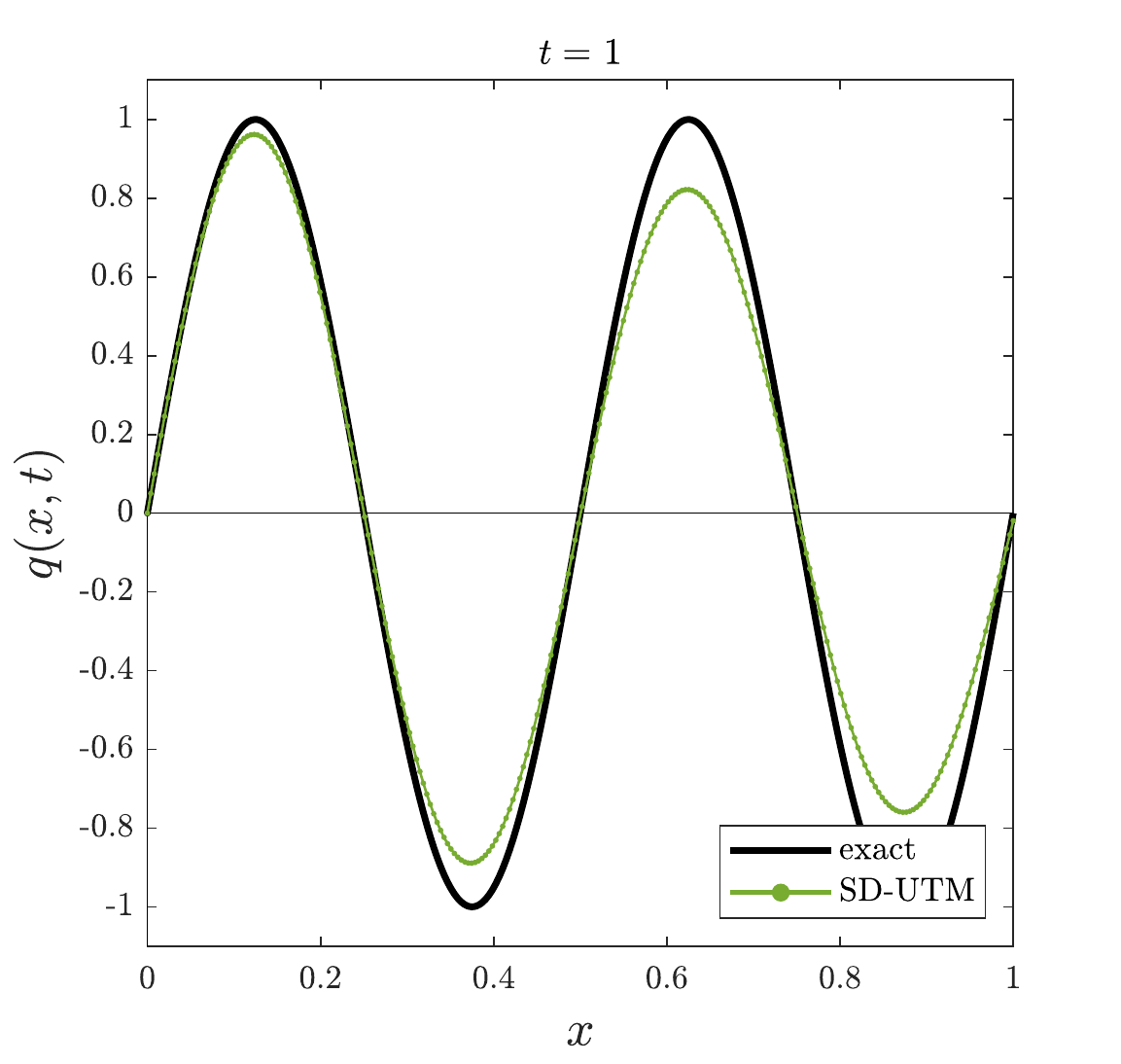}
	\end{subfigure}
	\caption{Several time slices for the solution to IBVP \eqref{advec2_numerical2_HL} with $h = 0.004$.}
	\label{advec2_numerical2_t}
\end{figure}

	As an example, we examine the IBVP
	\begin{equation}
\begin{dcases}
	q_t = - q_{x},& x > 0,\, t > 0, \\
	q(x,0) = \phi(x) =  e^{-x} \sin \left(4 \pi x \right),& x > 0,  \\
	q(0,t) = u(t) = - \sin\left(4 \pi t\right),& t > 0,
\end{dcases}
\label{advec2_numerical2_HL}
\end{equation}	
	where the continuous solution is given by 
	\[
	q(x,T) = \begin{dcases}
		u^{(0)}\left(T-x\right)\, , \quad &0<x<T, \\
		\phi(x-T)\,, \quad & x > T.
	\end{dcases}
	\]
	Applying the semi-discrete solution \eqref{soln_advec2_backward} gives Figure \ref{advec2_UTM2_HL}, similar to Figure \ref{advec1_UTM1_HL} for this IBVP, illustrating the qualitative behavior of the advection equation and the expected $\mathcal{O}(h)$ error as $h \rightarrow 0$.
	Since we are using a purely-one sided stencil, the standard backward stencil, we expect the solution to be better approximated near $x = 0$ and less so for larger $x$. From the stencil \eqref{advec2_backward}, we find the convection-diffusion PDE $p_t = - c\, p_x + (c h/2) p_{xx}$ as its modified PDE. Like \eqref{advec_forward_modified_eqn}, this modified equation is approximately solved by the semi-discrete solution with second-order accuracy. The presence of the dissipative term implies \eqref{soln_advec2_backward} advects the initial and boundary data at the appropriate speed, but with $\mathcal{O}(h)$ damping as time progresses. Indeed, Figure \ref{advec2_numerical2_t} displays the dissipation present in the stencil, manifested in its modified PDE, away from the boundary with a series of plots for various $t$ and $h = 0.004$. 
\newpage


\subsection{Centered Discretization of $\bms{q_t = -c\,q_{x}}$}\label{advec2_centered_halfline}
\begin{figure}[tb]
		\begin{center}
			\includegraphics[width=0.4\linewidth]{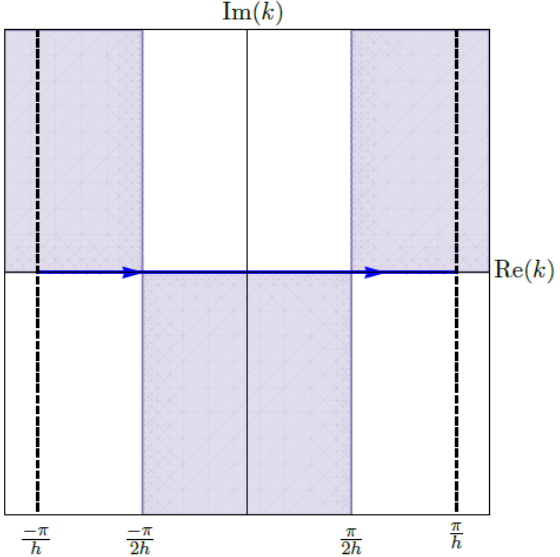}			
			\caption{The shaded regions depict where $\text{Re}(-W) \leq 0$ and $e^{-WT}$ is bounded with the dispersion relation \eqref{W_advec2_centered}.}
			\label{advec_cent_W}
		\end{center}
	\end{figure}

	We consider the same problem as in Section \ref{advec2_backward_halfline}, but using the standard centered discretization:
	\begin{align}
		\dot{q}_n(t) = - c\,\left(\frac{q_{n+1}(t) - q_{n-1}(t)}{2h}\right).
		\label{advec2_centered}
	\end{align}
	With slightly more work, the local relation is 
	\begin{align}
		\partial_t \left(e^{-iknh} e^{Wt} q_n \right) &= \frac{-c}{2h}\Delta \left( e^{-iknh} e^{Wt} q_{n-1} + e^{-ik(n-1)h} e^{Wt} q_{n}\right),
		\label{LR_advec2_centered}
	\end{align}
	with dispersion relation 
	\begin{equation}
		W(k) = c\,\frac{e^{ikh} - e^{-ikh}}{2h}= \frac{- c\,\sin(kh)}{i h}.
		\label{W_advec2_centered}
	\end{equation}
	In this case, the dispersion relation has the trivial symmetry $\nu_0 = k$ and one nontrivial symmetry
	$$\nu_1(k) = - k - \frac{\pi}{h},$$
	up to periodic copies. Since we have information at $q(0,t) \equiv q_0(t)$, we take an infinite sum starting at $n = 1$ and a time transform to obtain the global relation as
	\begin{align}
		e^{WT} \hat{q}(k,T) - \hat{q}(k,0) - c\,\left[ \frac{e^{-i k h} f_0 + f_1}{2} \right] &= 0,\quad \text{Im}(k) \leq 0.
		\label{GR_advec2_centered}
	\end{align}
	Taking the inverse transform, we obtain the ``solution''
	\begin{align}\begin{split}
		q_n(T) &= \frac{1}{2\pi} \int_{-\pi/h}^{\pi/h} e^{iknh} e^{-WT} \hat{q}(k,0)\,dk + \frac{c}{2\pi} \int_{-\pi/h}^{\pi/h} e^{iknh} e^{-WT} \left[ \frac{e^{-i k h} f_0 + f_1}{2} \right]\,dk.
		\label{soln1_advec2_centered}
	\end{split}\end{align}
	``Solution'' \eqref{soln1_advec2_centered} contains the unknown $f_{1}(W,T)$, but as seen in Figure \ref{advec_cent_W}, we cannot deform off the real line to remove this dependence as we did with one-sided discretization stencils.
	Nonetheless, the global relation \eqref{GR_advec2_centered} with $k \rightarrow \nu_1(k)$ is valid for $\text{Im}\left(\nu_1\right) \leq 0$, \textit{i.e.}, $\text{Im}(k) \geq 0$, and can be used to remove the unknown without the need to deform. Replacing $k \, \rightarrow \, \nu_1$ in the global relation \eqref{GR_advec2_centered} and substituting
	$$f_{-1} = - e^{-i \nu_1 h} f_0 + \frac{2}{c} \left[ e^{WT} \hat{q}\left(\nu_1,T\right) - \hat{q}\left(\nu_1,0\right) \right],\quad \text{Im}(k) \geq 0,$$
	in the ``solution'' \eqref{soln1_advec2_centered}, we find
	\begin{align}
		\begin{split}
		q_n(T) &= \frac{1}{2\pi} \int_{-\pi/h}^{\pi/h} e^{iknh} e^{-WT} \hat{q}(k,0)\,dk + \frac{1}{2\pi} \int_{-\pi/h}^{\pi/h} e^{iknh} e^{-WT} \left[c\, \cos\left( kh \right)f_0 -\hat{q}\left(\nu_1,0\right) \right] \,dk \\
		&\quad\, + \frac{1}{2\pi} \int_{-\pi/h}^{\pi/h} e^{iknh} \hat{q}\left(\nu_1,T\right) \,dk,
		\label{soln2_advec2_centered}
		\end{split}
	\end{align}
	after simplification.	 
	
	Removing one unknown from the ``solution,'' we have introduced another, $\hat{q}\left(\nu_1,T\right)$, a transform of the solution at time $T$. It is crucial to point out that this last integral term does not have the exponential factor $e^{-WT}$. We can eliminate $\hat{q}\left(\nu_1,T\right)$ from our ``solution'' as in the continuous UTM, or determine its contribution if it is nonzero. To do so, we substitute the definition of the transform into the integral term:
	\begin{align*}
		\frac{1}{2\pi} \int_{-\pi/h}^{\pi/h} e^{iknh} \hat{q}\left(\nu_1,T\right) \,dk &= \frac{1}{2\pi} \int_{-\pi/h}^{\pi/h} e^{iknh} \left[h \sum_{m=1}^\infty e^{-i \nu_1 mh} q_m(T) \right] \,dk = \sum_{m=1}^\infty (-1)^m C(n+m) q_m(T),
	\end{align*} 
	where 
	$$C(n) = \frac{h}{2\pi} \int_{-\pi/h}^{\pi/h} e^{iknh} \,dk.$$
	Applying the first approach from Subsection \ref{advec_forward_halfline} implies that for $n > 0$, $C(n) = 0$ via periodicity, 
	$$\frac{1}{2\pi} \int_{-\pi/h}^{\pi/h} e^{iknh} \hat{q}\left(\nu_1,T\right) \,dk  = 0.$$
	Therefore, the solution to the centered-discretized advection equation $q_t = -c\, q_x$ on the half-line with a Dirichlet boundary condition is
	\begin{align}
	\begin{split}
		q_n(T) &= \frac{1}{2\pi} \int_{-\pi/h}^{\pi/h} e^{iknh} e^{-WT} \hat{q}(k,0)\,dk + \frac{1}{2\pi} \int_{-\pi/h}^{\pi/h} e^{iknh} e^{-WT} \left[c\, \cos\left( kh \right)f_0 -\hat{q}\left(\nu_1,0\right) \right] \,dk .
		\label{soln_advec2_centered}
		\end{split}
	\end{align}
	
	From the stencil \eqref{advec2_centered}, we find the modified PDE $p_t =  - c \,p_x - (ch^2/6) p_{xxx},$ preserving the correct advection speed, but including a dispersive term. Hence, the solution profiles disperse as time progresses. Since the dispersion coefficient is $\mathcal{O}\left(h^2\right)$, these effects are minimal for practical $h \ll 1$.  
	
	The continuum limit of this semi-discrete solution is less straightforward than the continuum limit of the backward-discretized solution \eqref{soln_advec2_backward}. The dispersion relation converges to the continuous one: $\lim_{h \rightarrow 0}W(k) = \tilde{W}(k)$. It is clear that the first integral term in the semi-discrete solution \eqref{soln_advec2_centered} converges to the first integral term in the continuous solution \eqref{soln_advec2_cont}. Since $\lim_{h \rightarrow 0}\cos(kh) = 1$, the boundary component in the discrete solution also converges to its continuous counterpart. 
	
	The solution to the continuous problem requires no additional symmetries, so we expect the integral term containing $\hat{q}\left(\nu_1,0\right)$ to vanish. Note that
	\begin{align*}
		\hspace{-15pt}\hat{q}\left(\nu_1,t\right) &= h \sum_{m=1}^\infty e^{-i \nu_1 mh} q_m(t) = h \sum_{m=1}^\infty (-1)^m e^{i k mh} q_m(t) = - h \sum_{u=1}^\infty e^{2 i k u h} q_{2 u}(t) + h \sum_{v=0}^\infty e^{(2v + 1)i k h} q_{2v+1}(t),
	\end{align*}
	after separating the even and odd indexed terms. Similar to $x_n = n h$, we introduce dummy variables $w_u = u h$ and $y_v = v h$, so that 
	\begin{align*}
		\lim_{h \rightarrow 0} \hat{q}\left(\nu_1,t\right) &= - \lim_{h \rightarrow 0} h \sum_{u=1}^\infty e^{2 i k u h} q_{2 u}(t) + \lim_{h \rightarrow 0}  h \sum_{v=0}^\infty e^{2 i k v h}e^{ikh} q_{2v+1}(t) \\
		&= - \int_{0}^\infty e^{2 i k w} q(2w,t)\,dw + \int_{0}^\infty e^{2 i k y} q(2y,t)\,dy = 0.
	\end{align*}
	Thus,
	$$ \lim_{h \rightarrow 0}\, \frac{-1}{2\pi} \int_{-\pi/h}^{\pi/h} e^{iknh} e^{-WT} \hat{q}\left(\nu_1,0\right) \,dk = 0,$$
	and we recover the continuous solution \eqref{soln_advec2_cont}.


	\subsection{Higher-Order One-Sided Discretization of $\bms{q_t = -c\,q_{x}}$}\label{advec2_highorder_halfline}
	There exist higher-order discretizations that appropriately incorporate the nontrivial symmetries to remove unknowns, where the steps in the semi-discrete UTM become more intricate and tedious, yet remain systematic. In some cases, however, the nontrivial symmetries are unusable, but a solution can still be obtained. Consider the second-order discretized advection equation
	\begin{equation}
		\dot{q}_n(t) = -c\left(\frac{q_{n-2}(t) - 4 q_{n-1}(t) + 3q_{n}(t)}{2h}\right).
		\label{advec2_backward2}
	\end{equation}
	Following the usual steps, the global relation is
	\begin{align}
		 e^{WT} \hat{q}(k,T) - \hat{q}(k,0) - c \left( \frac{ 4 e^{-ikh} f_{0} - e^{-2ikh} f_{0} - e^{-ikh} f_{-1} }{2}\right)&= 0, \quad \text{Im}(k) \leq 0,
		\label{GR_advec2_backward2}
	\end{align}
	with dispersion relation
	\begin{equation}
		W(k) = c \,\frac{e^{-2ikh} - 4 e^{-ikh} + 3}{2h},
		\label{W_advec2_backward2}
	\end{equation}
	and nontrivial symmetry 
	$$\nu_1(k) = \frac{i}{h} \ln\left(4 - e^{-ikh}\right).$$
	Solving for $\hat{q}(k,T)$ and taking the inverse transform, we obtain
	\begin{align}\begin{split}
		q_n(T) &= \frac{1}{2\pi} \int_{-\pi/h}^{\pi/h} e^{iknh} e^{-WT} \hat{q}(k,0)\,dk \\
		&\quad\, + \frac{c}{2\pi} \int_{-\pi/h}^{\pi/h} e^{iknh} e^{-WT} \left( \frac{ 4 e^{-ikh} f_{0} - e^{-2ikh} f_{0} - e^{-ikh} f_{-1} }{2}\right)\,dk.
		\label{soln1_advec2_backward2}
	\end{split}\end{align}
	\begin{figure}[tb]
		\begin{center}
			\includegraphics[width=0.55\linewidth]{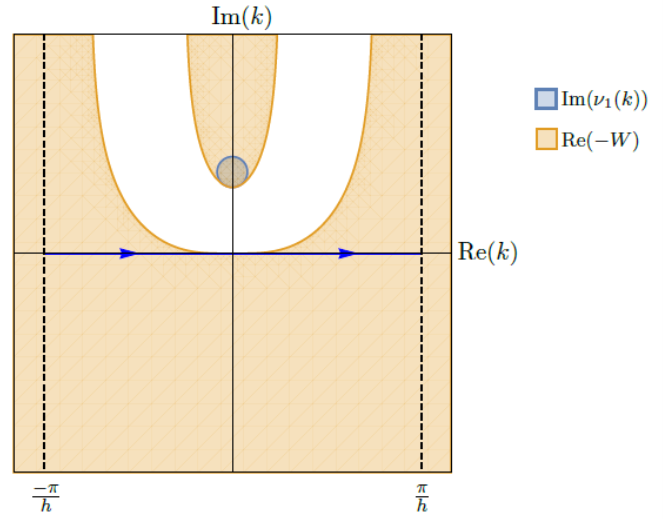}
			\caption{The blue shaded region denotes where the global relation \eqref{GR_advec2_backward2} is valid with $k \rightarrow \nu_1(k)$ and the orange shaded regions depict where $\text{Re}(-W) \leq 0$ and $e^{-WT}$ is bounded with the dispersion relation \eqref{W_advec2_backward2}.}
			\label{symm_advec2_backward2}
		\end{center}
	\end{figure}
	With the given Dirichlet boundary condition, we encounter the ghost point $f_{-1}(W,T)$, which is unknown. The global relation \eqref{GR_advec2_backward2} is valid in the blue shaded region of Figure \ref{symm_advec2_backward2} with $k \rightarrow \nu_1(k)$, but since the path on the real line cannot be deformed to this region, the nontrivial symmetry cannot be used to eliminate the unknown $f_{-1}(W,T)$. Instead, we return to the continuous problem \eqref{advec2_prob}, where the PDE itself gives the Neumann boundary condition from the Dirichlet condition:
	\begin{align}
		q_x(0,t) = \frac{-1}{c}\,q_{t}(0,t) = \frac{-1}{c}\,\frac{\partial}{\partial t} q(0,t) = \frac{-u'(t)}{c} = \frac{-v(t)}{c}, \quad v(t) = u'(t).
		\label{advec2_neumann_bc}
	\end{align}
	In order to remove $f_{-1}(W,T)$ without introducing new unknowns, we discretize the Neumann condition using the standard backward stencil:
	\begin{align}
		&&\hspace{-75pt}\frac{q_{0}(t) - q_{-1}(t)}{h} &= \frac{-v(t)}{c} \notag \\ 
		\Rightarrow &&\hspace{-75pt} \frac{f_{0} - f_{-1}}{h} &= \frac{-V}{c}, \quad\quad\quad V(W,T) = \int_0^T e^{Wt} v(t)\,dt \notag\\ 
		\Rightarrow &&\hspace{-75pt} f_{-1} &= f_{0} + \frac{h}{c} V.
		\label{advec2_neumann_bc_f}
	\end{align}
	Substituting in \eqref{soln1_advec2_backward2}, we obtain the solution
	\begin{align}\begin{split}
		q_n(T) &= \frac{1}{2\pi} \int_{-\pi/h}^{\pi/h} e^{iknh} e^{-WT} \hat{q}(k,0)\,dk \\
		&\quad\, + \frac{c}{2\pi} \int_{-\pi/h}^{\pi/h} e^{iknh} e^{-WT} \left[ \left(\frac{ 3 e^{-ikh} - e^{-2ikh}}{2}\right) f_{0} - \left(\frac{h e^{-ikh}}{2c}\right)V\right]\,dk.
		\label{soln_advec2_backward2}
	\end{split}\end{align}
	Of course, the Neumann boundary condition was discretized to first-order accuracy, so the accuracy of \eqref{soln_advec2_backward2} is $\mathcal{O}(h)$. As before, as $h \rightarrow 0$, the semi-discrete solution \eqref{soln_advec2_backward2} converges to \eqref{soln_advec2_cont}, and the semi-discrete solution correctly loses dependence on the Neumann boundary condition in the continuum limit.
	
	\begin{remark}
	Let us reconsider the backward discretization \eqref{advec2_backward} with dispersion relation \eqref{W_advec2_backward} and no nontrivial symmetries. With the forward transform $\hat{q}(k,t)$ starting at $n = 0$ instead of $n = 1$, we derive the global relation
	\begin{equation}
		e^{WT} \hat{q}(k,T) - \hat{q}(k,0) - c f_{-1} = 0, \quad \text{Im}(k) \leq 0,
		\label{GR_advec2_backward_transform_n0}
	\end{equation}
	with ``solution'' 
	\begin{align}
		q_n(T) &= \frac{1}{2\pi} \int_{-\pi/h}^{\pi/h} e^{iknh} e^{-WT} \hat{q}(k,0)\,dk + \frac{c}{2\pi} \int_{-\pi/h}^{\pi/h} e^{iknh} e^{-WT} f_{-1}\,dk,
		\label{soln_advec2_backward_n0}
	\end{align}
	depending on $q_{-1}(t)$, which is not directly provided by the IBVP \eqref{advec2_prob}. As above, the advection equation itself gives the Neumann boundary condition \eqref{advec2_neumann_bc} from the given Dirichlet condition and allows us to remove the dependence on $f_{-1}(W,T)$ in \ref{soln_advec2_backward_n0}. Substituting \eqref{advec2_neumann_bc_f}, the solution is
	\begin{align}\begin{split}
		q_n(T) &= \frac{1}{2\pi} \int_{-\pi/h}^{\pi/h} e^{iknh} e^{-WT} \hat{q}(k,0)\,dk + \frac{c}{2\pi} \int_{-\pi/h}^{\pi/h} e^{iknh} e^{-WT} f_{0}\,dk + \frac{h}{2\pi} \int_{-\pi/h}^{\pi/h} e^{iknh} e^{-WT} V \,dk.
		\label{soln_advec2_backward_n0_final}
	\end{split}\end{align}
	In the continuum limit, the last integral term vanishes and we recover \eqref{soln_advec2_cont}. Both \eqref{soln_advec2_backward} and \eqref{soln_advec2_backward_n0_final} are solutions to the backward-discretized advection equation \eqref{advec2_backward} with $q_0(t)$ data, but the transforms $\hat{q}(k,t)$ are defined differently by a shift in the starting index. Using the global relations \eqref{GR_advec2_backward_transform} and \eqref{GR_advec2_backward_transform_n0}, one can show that the solutions \eqref{soln_advec2_backward} and \eqref{soln_advec2_backward_n0_final} are equal.
	
	\end{remark}


\section{The Heat Equation}

\subsection{Centered Discretization of $\bms{q_t = q_{xx}}$ with Dirichlet boundary condition}\label{heat_centered_halfline_sec}
	Consider the problem
		\begin{equation}\begin{dcases}
		q_t = q_{xx},& x > 0,\, t > 0, \\
		q(x,0) = \phi(x),& x > 0,\\
		q(0,t) = u(t),& t > 0,
		\label{heat_prob}
	\end{dcases}\end{equation}
	with one Dirichlet boundary condition. We write the centered-discretized heat equation as
	\begin{equation}
		\dot{q}_n(t) = \frac{q_{n+1}(t) - 2 q_n(t) + q_{n-1}(t)}{h^2}.
		\label{heat_centered}
	\end{equation}
	Carrying out similar steps as before, the local relation is
	\begin{align}
		\partial_t \left(e^{-iknh} e^{Wt} q_n \right) &= \frac{1}{h^2}\Delta \left(e^{-ik(n-1)h} e^{Wt} q_{n} - e^{-iknh} e^{Wt} q_{n-1} \right),
		\label{LR_heat_centered}
	\end{align}
	where
	\begin{equation}
		W(k) = \frac{2 - e^{ikh} - e^{-ikh}}{h^2} = \frac{2}{h^2} \left[ 1 - \cos(kh) \right],
		\label{W_heat_centered}
	\end{equation}
	with the nontrivial symmetry $\nu_1(k) = -k$ up to periodic copies. We find the global relation by summing from $n = 1$ and integrating in time:
	\begin{align}
		e^{WT} \hat{q}(k,T) - \hat{q}(k,0) - \left[ \frac{e^{-ikh} f_{0} - f_1}{h} \right] &= 0, \quad \text{Im}(k) \leq 0.
		\label{GR_heat_centered}
	\end{align}
	Inverting, we obtain the ``solution'' formula
	\begin{align}
		q_n(T) &= \frac{1}{2 \pi} \int_{-\pi/h}^{\pi/h} e^{iknh} e^{-WT}\hat{q}(k,0)\,dk + \frac{1}{2 \pi} \int_{-\pi/h}^{\pi/h} e^{iknh} e^{-WT}\left[ \frac{e^{-ikh} f_{0} - f_1}{h} \right]\,dk,
		\label{soln1_heat_centered}
	\end{align}
	which depends on the unknown $f_1(W,T)$. Using $\nu_1(k)$ in the global relation \eqref{GR_heat_centered} gives
	\begin{align}
		f_1 &= e^{ikh} f_{0} - h\left[ e^{WT} \hat{q}(-k,T) - \hat{q}(-k,0) \right], \quad \text{Im}(k) \geq 0.
		\label{GR_heat_centered_1}
	\end{align}
	We substitute \eqref{GR_heat_centered_1} into \eqref{soln1_heat_centered}, so that
	\begin{align}
		\begin{split}
			q_n(T) &= \frac{1}{2 \pi} \int_{-\pi/h}^{\pi/h} e^{iknh} e^{-WT}\hat{q}(k,0)\,dk - \frac{1}{2 \pi} \int_{-\pi/h}^{\pi/h} e^{iknh} e^{-WT}\left[ \hat{q}(-k,0) + \frac{ 2i \sin(kh) }{h}f_{0} \right]\,dk \\
		&\quad\,+ \frac{1}{2 \pi} \int_{-\pi/h}^{\pi/h} e^{iknh} \hat{q}(-k,T) \,dk. 
		\label{soln2_heat_centered} 
		\end{split} 
	\end{align}
	
		Removing the boundary term $f_{1}(W,T)$, we have introduced the transform of the solution at $t = T$ in the third integral of \eqref{soln2_heat_centered}. Using the definition of the transform, 
	\begin{align*}
		\frac{1}{2 \pi} \int_{-\pi/h}^{\pi/h} e^{iknh} \hat{q}(-k,T) \,dk = \frac{1}{2 \pi} \int_{-\pi/h}^{\pi/h} e^{iknh} \left[h \sum_{m=1}^\infty e^{ikmh} q_m(T)\right] \,dk = \sum_{m=1}^\infty q_m(T) \left[\frac{h}{2 \pi} \int_{-\pi/h}^{\pi/h} e^{ik(n+m)h}\,dk\right].
	\end{align*}
	For $n > 0$, the integral vanishes due to periodicity, so that 
	$$\frac{1}{2 \pi} \int_{-\pi/h}^{\pi/h} e^{iknh} \hat{q}(-k,T) \,dk =0,$$
	and the solution to \eqref{heat_centered} with the Dirichlet boundary condition is written as
	\begin{equation}
	\begin{aligned}
		q_n(T) &= \frac{1}{2 \pi} \int_{-\pi/h}^{\pi/h} e^{iknh} e^{-WT}\hat{q}(k,0)\,dk  - \frac{1}{2 \pi} \int_{-\pi/h}^{\pi/h} e^{iknh} e^{-WT}\left[ \hat{q}(-k,0) + \frac{2i \sin(kh) f_{0}}{h} \right]\,dk.
		\label{soln_heat_centered}
	\end{aligned}
	\end{equation}
	
	Solving the IBVP \eqref{heat_prob} via the continuous UTM gives the dispersion relation $\tilde{W}(k) = k^2$ with nontrivial symmetry $\tilde{\nu}_{1}(k) = -k$ \cite{bernard_fokas}. The solution is
	 \begin{align}
		q(x,T) &= \frac{1}{2 \pi} \int_{-\infty}^{\infty} e^{ikx} e^{-\tilde{W}T} \hat{q}(k,0)\, dk - \frac{1}{2 \pi} \int_{-\infty}^{\infty} e^{ikx} e^{-\tilde{W}T} \left[\hat{q}(-k,0)  + 2ik F_0 \right]\,dk.
		\label{soln_heat_cont}
	\end{align}
	Taking the continuum limit, \eqref{soln_heat_centered} converges to \eqref{soln_heat_cont}, since $\lim_{h \rightarrow 0} W(k) = k^2 = \tilde{W}(k)$.
	
\begin{figure}[tb]	
		\raggedleft
		\begin{subfigure}[t]{.45\textwidth}
			\centering
  			\includegraphics[width=0.95\linewidth]{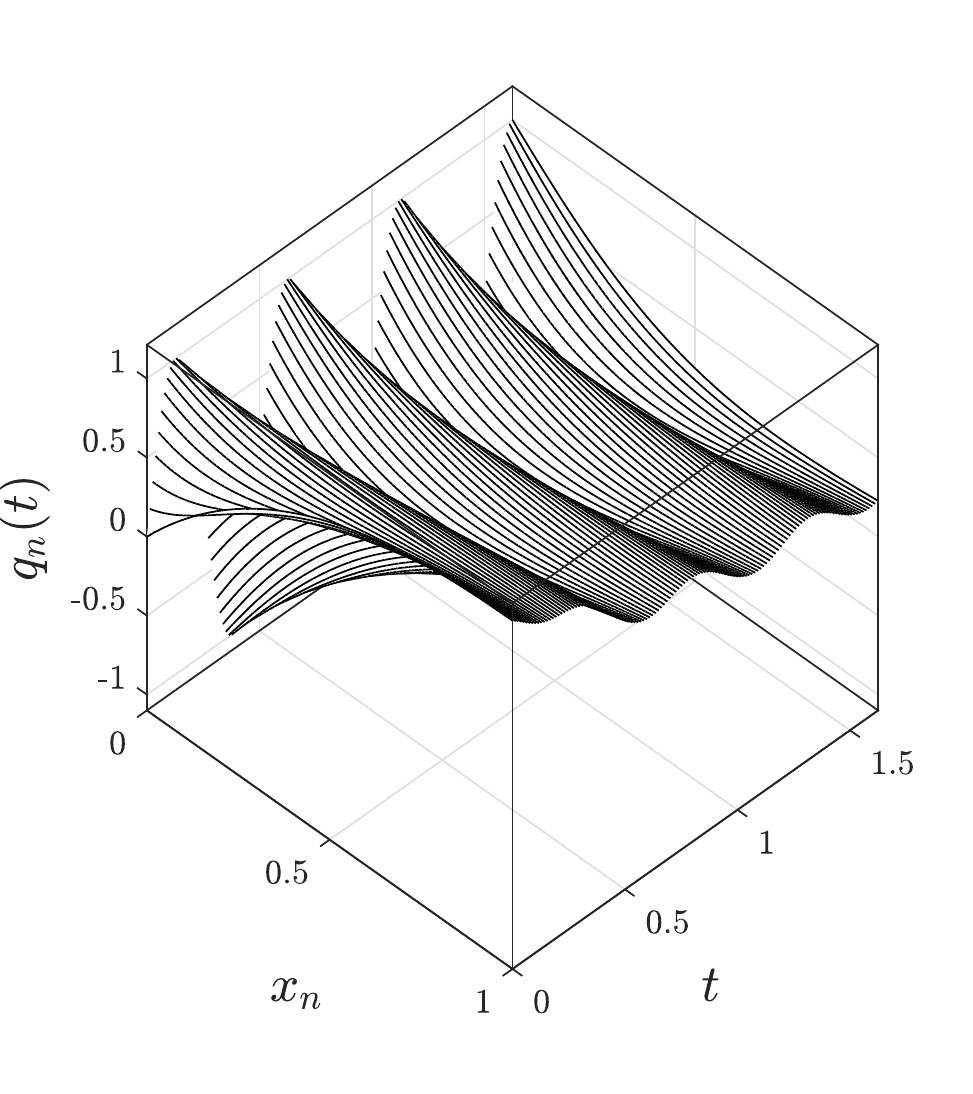}
  			\caption{}
  			\label{}
		\end{subfigure}\hfill 
		\begin{subfigure}[t]{.45\textwidth}
			\centering
  			\includegraphics[width=1\linewidth]{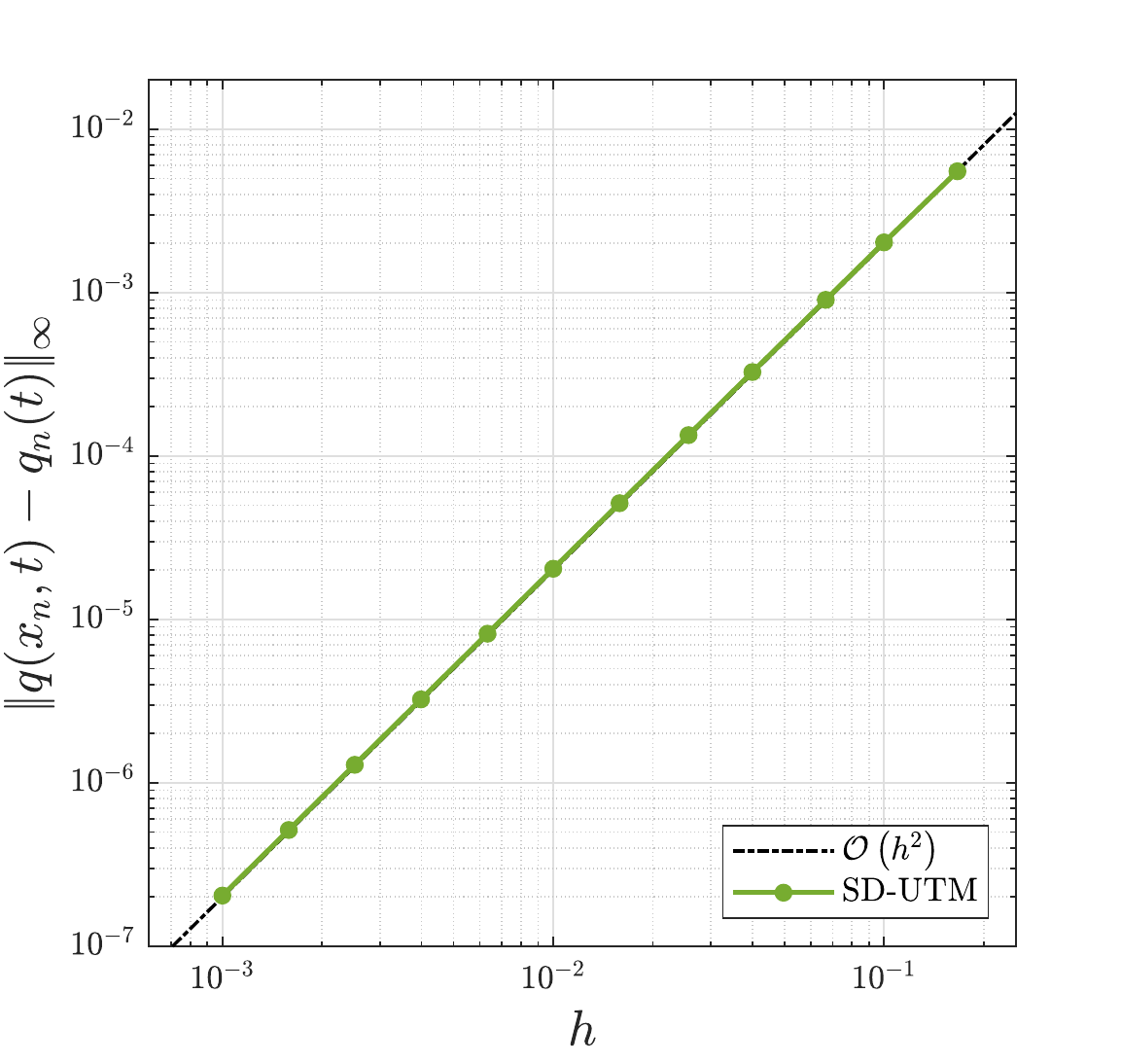}
  			\caption{}
  			\label{}
		\end{subfigure}	
		\caption{(a) The semi-discrete solution \eqref{soln_heat_centered} evaluated at various $t$ with $h = 0.01$. (b) Error plot of the semi-discrete solution \eqref{soln_heat_centered} relative to the exact solution as $h \rightarrow 0$ with $t = 1.625$.}
		\label{heat_UTM2_HL}
	\end{figure}

	As an example, the solution to the IBVP
\begin{equation}
\begin{dcases}
	q_t = q_{xx},& x > 0,\, t > 0, \\
	q(x,0) = \phi(x) = 3 x e^{-x},& x > 0,  \\
	q(0,t) = u(t) = \sin\left(4 \pi t\right),& t > 0.
\end{dcases}
\label{heat_numerical2_HL}
\end{equation}	
	is written in terms of error functions. Deriving the modified PDE from the centered stencil \eqref{heat_centered}, we find that solution \eqref{soln_heat_centered} is a fourth-order accurate approximation to the solution of the dissipative PDE 
	\begin{equation}
		p_t = p_{xx} +\frac{h^2}{12} p_{4x}.
		\label{heat_centered_modified_eqn}
	\end{equation} 
	The presence of the higher-order dissipation term $p_{4x}$ causes high-frequency oscillations to be damped for any $t >0$. The original heat equation is also dissipative, but solution \eqref{soln_heat_centered} might overdamp in scenarios where the initial data contains high-frequency oscillations or the boundary condition oscillates in time with large amplitude. Although the dissipation coefficient of $p_{4x}$ is $\mathcal{O}\left(h^2\right)$, the overdamping nature can be troublesome for a practical $h \ll 1$ as $t$ increases, but this can be counteracted by decreasing $h$. With the SD-UTM solution \eqref{soln_heat_centered}, the left plot of Figure \ref{heat_UTM2_HL} shows the gradual decay away from the $x=0$ and $t = 0$ boundaries as time increases and the right plot shows the expected $\mathcal{O}\left(h^2\right)$ error as $h \rightarrow 0$.


\subsection{Centered Discretization of $\bms{q_t = q_{xx}}$ with Neumann boundary condition} \label{neumann_halfline}
	
	We consider the continuous half-line problem:
		\begin{equation}\begin{dcases}
		q_t = q_{xx},& x > 0,\, t > 0, \\
		q(x,0) = \phi(t),& x > 0,\\
		q_x(0,t) = v(t),& t > 0,
		\label{heat_prob_N}
	\end{dcases}\end{equation}
	with a Neumann boundary condition. How do we discretize this condition so that we may employ it with the centered-discretized equation \eqref{heat_centered}? This choice often leads to instabilities in finite-difference schemes, especially when dealing with higher-order problems \cite{cheema_thesis,randy,trefethen}. We show that the SD-UTM determines which discretizations we can choose.
	
	We proceed with the centered discretization \eqref{heat_centered} for the heat equation. This implies we retain the local relation \eqref{LR_heat_centered} and dispersion relation \eqref{W_heat_centered} with nontrivial symmetry $\nu_1(k) = -k$. We cannot use the global relation \eqref{GR_heat_centered}, because we assumed Dirichlet boundary data to obtain it. Now, we do not have information at $n = 0$ and we define our forward transform to start at $n = 0$:
	$$\hat{q}\left(k,t\right) =h \sum_{n=0}^{\infty} e^{-i k n h} q_{n}(t),$$
	directly affecting the global relation. From the local relation \eqref{LR_heat_centered}, 
	\begin{align}
		e^{WT} \hat{q}(k,T) - \hat{q}(k,0) - \left[ \frac{f_{-1} - e^{ikh} f_{0}}{h} \right] &= 0, \quad \text{Im}(k) \leq 0.
		\label{GR_heat_centered_N}
	\end{align}
	Solving for $\hat{q}(k,T)$ and inverting, we obtain
	\begin{align}
		q_n(T) = \frac{1}{2 \pi} \int_{-\pi/h}^{\pi/h} e^{iknh} e^{-WT}\hat{q}(k,0)\,dk + \frac{1}{2 \pi} \int_{-\pi/h}^{\pi/h} e^{iknh} e^{-WT}\left[ \frac{f_{-1} - e^{ikh} f_{0}}{h} \right]\,dk.
		\label{soln1_heat_centered_N}
	\end{align}
	
	The global relation \eqref{GR_heat_centered_N} and ``solution'' \eqref{soln1_heat_centered_N} depend on $q_0(t)$, as in the previous section, but also on $q_{-1}(t)$ instead of $q_1(t)$. For this Neumann IBVP, we do not have exact information at neither $q_{-1}(t)$ nor $q_{0}(t)$, and both $f_{-1}(W,T)$ and $f_0(W,T)$ are unknowns. Since we only have one nontrivial symmetry to remove one unknown, we must provide another equation involving $f_{-1}(W,T)$ or $f_0(W,T)$ in such a way as to not introduce any new unknowns. If we discretize the Neumann boundary condition $q_{x}(0,t) = v(t)$, the only approach is to use the standard backward stencil:
	$$\frac{q_{0}(t) - q_{-1}(t)}{h} = v(t).$$
	Upon taking a time transform, 
	\begin{equation}
		\frac{f_{0}(W,T) - f_{-1}(W,T)}{h} = V(W,T),\quad\quad  V(W,T) = \int_0^T e^{Wt} v(t)\,dt.
		\label{N_disctr}
	\end{equation}
	Of course, this discretization is $\mathcal{O}\left(h\right)$, while the centered discretization \eqref{heat_centered} for $q_{xx}$ is $\mathcal{O}\left(h^2\right)$. This suggests that the final semi-discrete solution will lose accuracy compared to the case with a Dirichlet boundary condition, but nonetheless converge to the continuous solution. The relation \eqref{N_disctr} becomes the second equation to remove the second unknown. Solving the system
	\[\begin{dcases*}
		e^{WT} \hat{q}(-k,T) - \hat{q}(-k,0) - \left[ \frac{f_{-1} - e^{-ikh} f_{0}}{h} \right] = 0, \\[6pt]
		\frac{f_{0} - f_{-1}}{h} =  V(t),
	\end{dcases*}\]
	for $f_{-1}(W,T)$ and $f_0(W,T)$ results in
	$$\frac{f_{-1} - e^{ikh} f_{0}}{h} =  e^{i k h} \hat{q}(-k,0) - \left(1+e^{i k h}\right) V(t) - e^{WT}e^{i k h} \hat{q}(-k,T), \quad \text{Im}(k) \geq 0.$$
	Since \eqref{soln1_heat_centered_N} has integration paths on the real line, direct substitution gives
	{\fontsize{11}{13.2} \selectfont\begin{align}
	\begin{split}
		\hspace{-25pt}q_n(T) &= \frac{1}{2 \pi} \int_{-\pi/h}^{\pi/h} e^{iknh} e^{-WT}\hat{q}(k,0)\,dk + \frac{1}{2 \pi} \int_{-\pi/h}^{\pi/h} e^{iknh} e^{-WT}\left[ e^{i k h} \hat{q}(-k,0) - \left(1+e^{i k h}\right) V(t) \right]\,dk \\
		\hspace{-25pt}&\quad\, - \frac{1}{2 \pi} \int_{-\pi/h}^{\pi/h} e^{ik(n+1)h} \hat{q}(-k,T) \,dk.
		\label{soln2_heat_centered_N}
	\end{split}
	\end{align}}
	As before, we introduced an unwanted term that depends on the transform of the solution. We show the contribution from this term is zero by substituting the definition for $\hat{q}(-k,T)$:
	\begin{align*}
		\frac{-1}{2 \pi} \int_{-\pi/h}^{\pi/h} e^{ik(n+1)h} \hat{q}(-k,T) \,dk &= \frac{-1}{2 \pi} \int_{-\pi/h}^{\pi/h} e^{ik(n+1)h} \left[h \sum_{m = 0}^{\infty} e^{i k m h} q_m(T) \right] \,dk \\
		&= \sum_{m = 0}^{\infty} \left[ \frac{-h}{2 \pi}  \int_{-\pi/h}^{\pi/h} e^{ik(n+1)h}  e^{i k m h}\,dk  \right] q_m(T)\\
		&= 0,
	\end{align*}
	where the integral vanishes by periodicity. The final solution to this IBVP with a Neumann boundary condition $v(t)$, discretized as above, is
	{\fontsize{11}{13.2} \selectfont\begin{align}
	\begin{split}
		\hspace{-25pt}q_n(T) &= \frac{1}{2 \pi} \int_{-\pi/h}^{\pi/h} e^{iknh} e^{-WT}\hat{q}(k,0)\,dk + \frac{1}{2 \pi} \int_{-\pi/h}^{\pi/h} e^{iknh} e^{-WT}\left[ e^{i k h} \hat{q}(-k,0) - \left(1+e^{i k h}\right) V(t) \right]\,dk.
		\label{soln_heat_centered_N}
	\end{split}
	\end{align}}
	
	Similarly as shown in \cite{bernard_fokas}, the solution representation for IBVP \eqref{heat_prob_N} using the continuous UTM is
	\begin{align}
		q(x,T) &= \frac{1}{2 \pi} \int_{-\infty}^{\infty} e^{ikx} e^{-\tilde{W}T} \hat{q}(k,0)\, dk + \frac{1}{2 \pi} \int_{-\infty}^{\infty} e^{ikx} e^{-\tilde{W}T} \big[ \hat{q}(-k,0) - 2 F_1 \big]\,dk.
		\label{soln_heat_cont_N}
	\end{align}
	Referencing \eqref{soln_heat_centered_N}, the continuum limits of the coefficients of $\hat{q}(-k,0)$ and $V(t)$ converge, where 
	$$\lim_{h\rightarrow 0} V(t) = \lim_{h \rightarrow 0} \int_0^T e^{Wt} v(t)\,dt = \int_0^T e^{\tilde{W}t} v(t)\,dt = F_1.$$
	
\begin{figure}[tb]	
		\raggedleft
		\begin{subfigure}[t]{.45\textwidth}
			\centering
  			\includegraphics[width=0.95\linewidth]{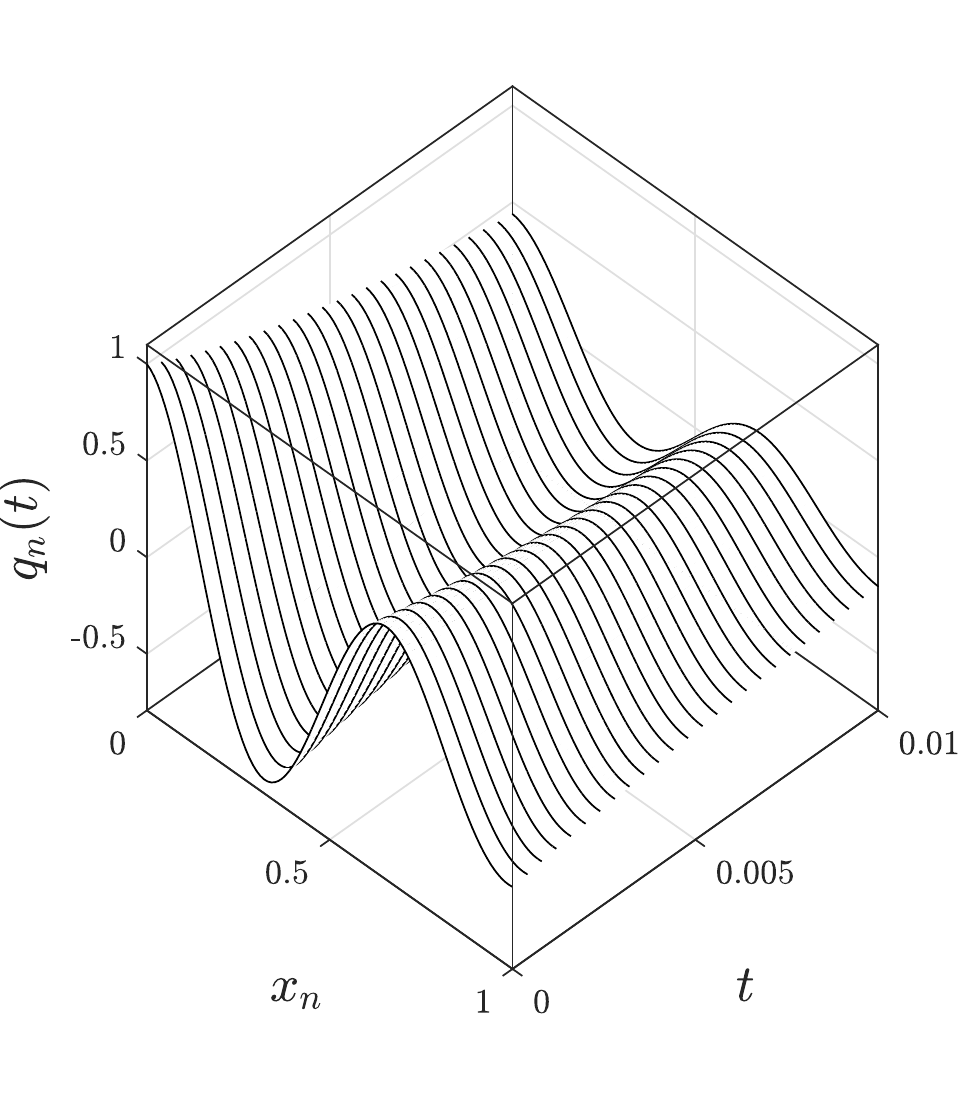}
  			\caption{}
  			\label{}
		\end{subfigure}\hfill 
		\begin{subfigure}[t]{.45\textwidth}
			\centering
  			\includegraphics[width=1\linewidth]{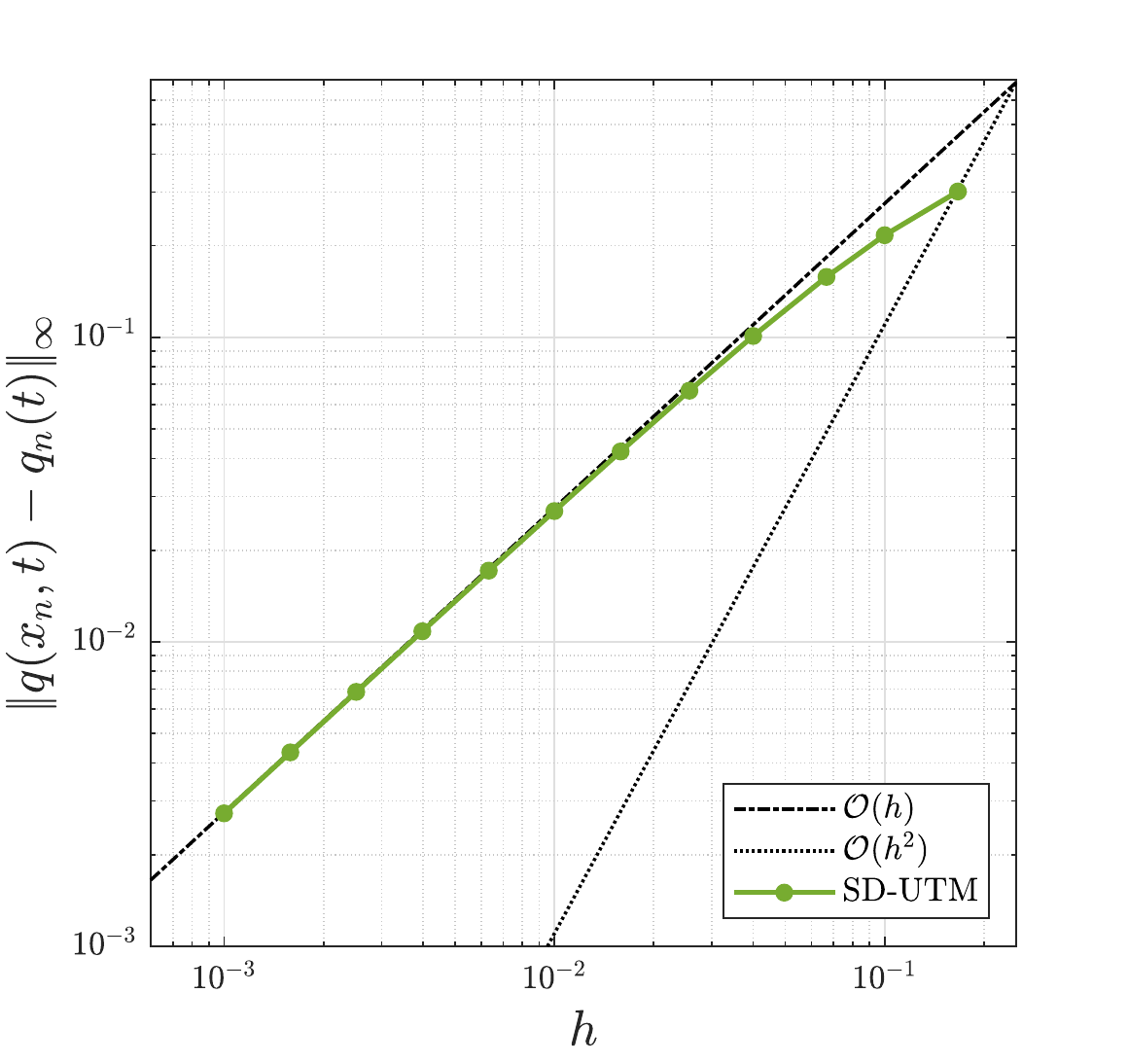}
  			\caption{}
  			\label{}
		\end{subfigure}	
		\caption{(a) The semi-discrete solution \eqref{soln_heat_centered_N} evaluated at various $t$ with $h = 0.01$. (b) Error plot of the semi-discrete solution \eqref{soln_heat_centered_N} relative to the exact solution as $h \rightarrow 0$ with $t = 0.01$.}
		\label{heat_N_UTM1_HL}
	\end{figure}

\begin{figure}[tb]
	\raggedleft
	\begin{subfigure}{.5\textwidth}
		\centering
  		\hspace*{-20pt}\includegraphics[width=0.9\linewidth]{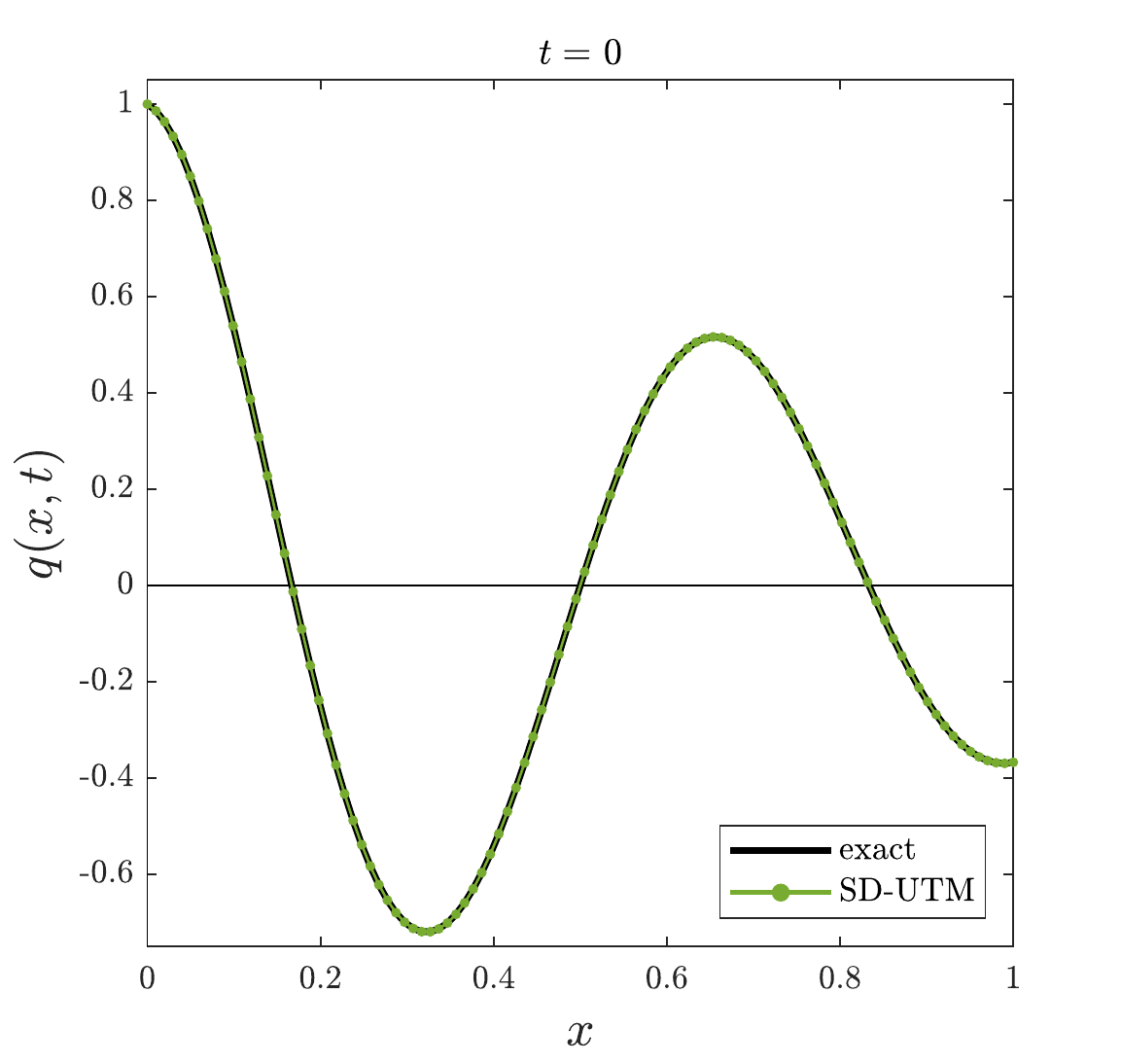}
	\end{subfigure}%
	\begin{subfigure}{.5\textwidth}
		\centering
  		\hspace*{0pt}\includegraphics[width=0.9\linewidth]{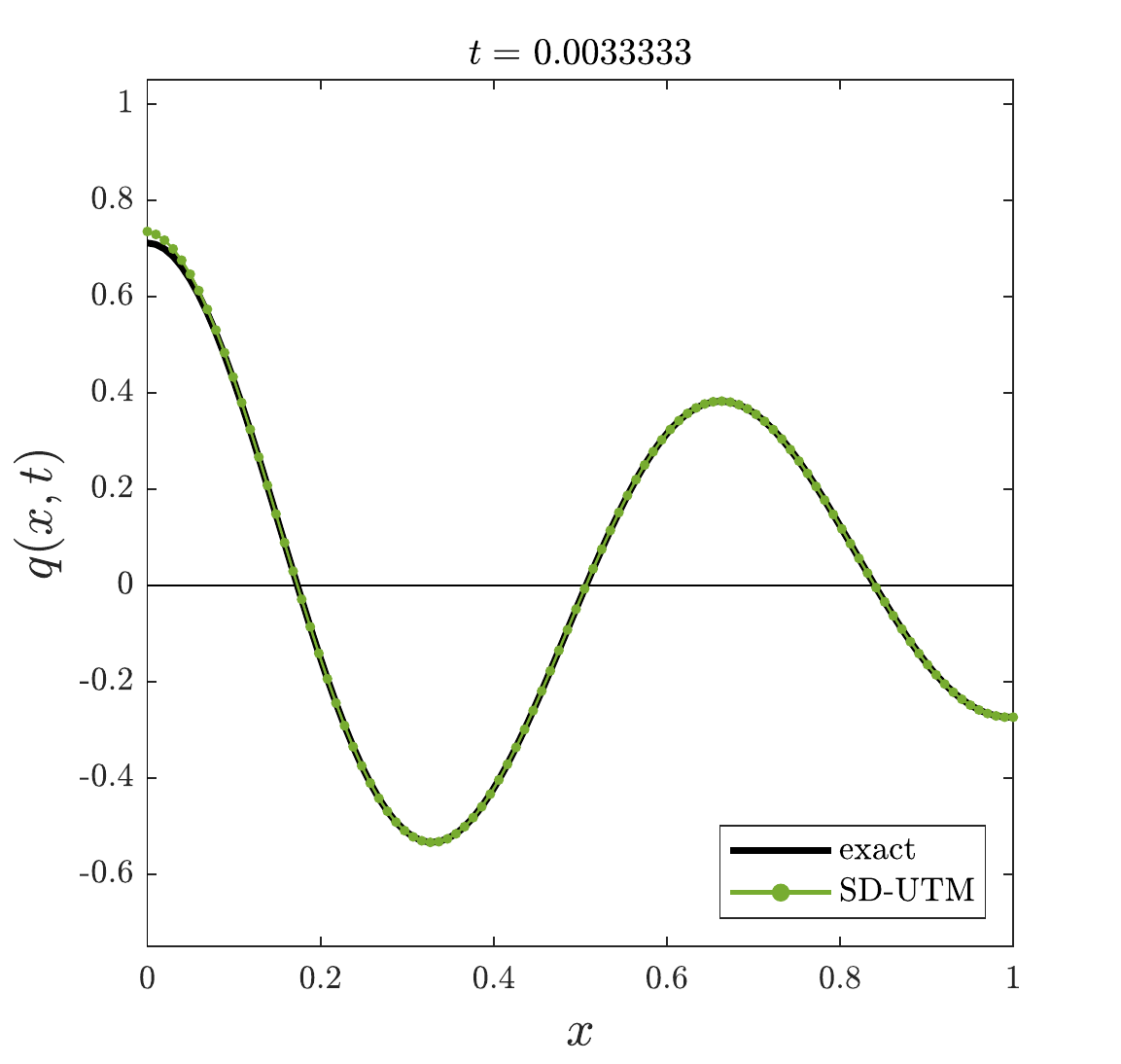}
	\end{subfigure}
	\begin{subfigure}{.5\textwidth}
		\centering
  		\hspace*{-20pt}\includegraphics[width=0.9\linewidth]{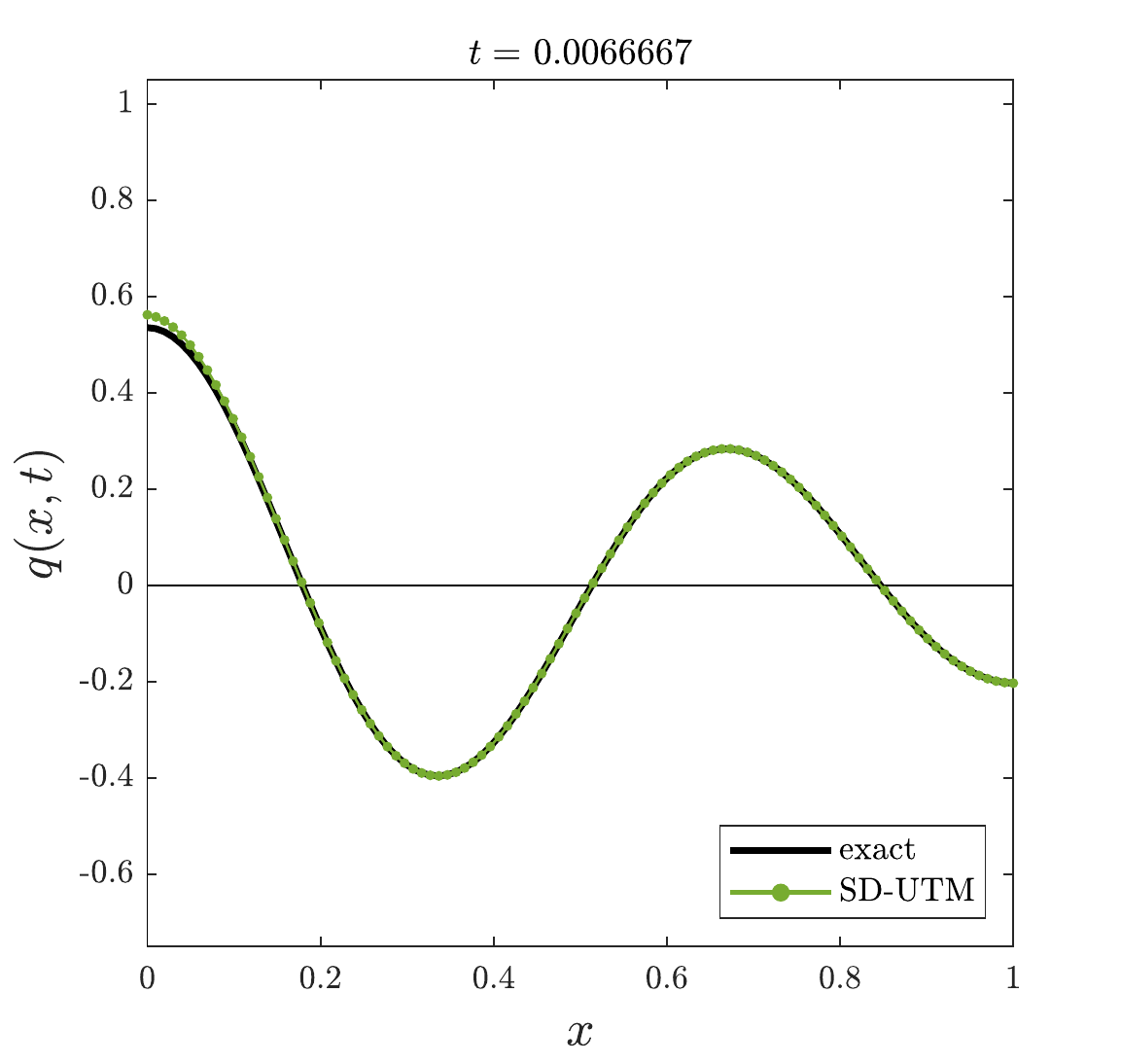}
	\end{subfigure}%
	\begin{subfigure}{.5\textwidth}
		\centering
  		\hspace*{0pt}\includegraphics[width=0.9\linewidth]{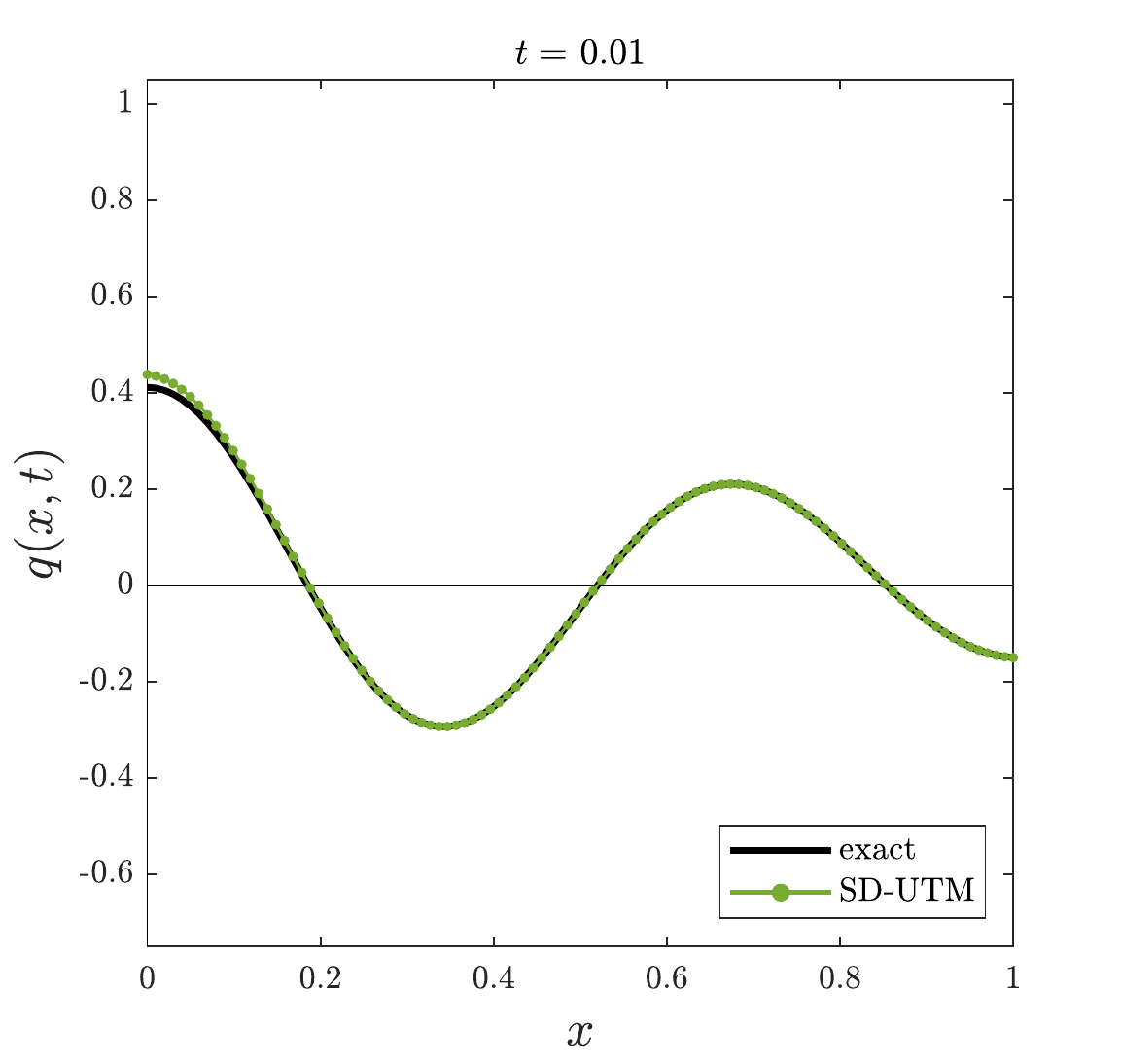}
	\end{subfigure}
	\caption{Several time slices for the solution to IBVP \eqref{heat_N_numerical1_HL} with $h = 0.01$.}
	\label{heat_N_numerical1_t}
\end{figure}

	As a concrete example, we examine the solution of the IBVP
\begin{equation}
\begin{dcases}
	q_t = q_{xx}, & x >0,\, t > 0, \\
	q(x,0) = \phi(x) = e^{-x} \cos\left(3 \pi x\right),& x > 0,  \\
	q_x(0,t) = v(t) = \tfrac{-1}{4\pi} \sin\left(4 \pi t\right),& t > 0.
\end{dcases}
\label{heat_N_numerical1_HL}
\end{equation}	
Again, the continuous solution is given in terms of error functions, while the semi-discrete solution is given by \eqref{soln_heat_centered_N}. The given Neumann data is discretized using the standard first-order accurate backward stencil, which reduces the overall accuracy of the solution from $\mathcal{O}(h^2)$ to $\mathcal{O}(h)$. Since the centered stencil \eqref{heat_centered} is used, solution \eqref{soln_heat_centered_N} is a fourth-order accurate approximation to the dissipative PDE \eqref{heat_centered_modified_eqn}. However, in general this is not the case, because of the discretization of the Neumann boundary condition. The modified equation for this backward discretization at $x = 0$ is $q_{x}(0,t) = v(t) - (h/2)q_{xx}(0,t).$ Unless the next several higher-order derivatives of the solution at the boundary are zero, the standard backward discretization we employed on $v(t)$ is $\mathcal{O}(h)$, and so is solution \eqref{soln_heat_centered_N}. This modified PDE at the boundary implies the loss of accuracy is visible in the solution profiles of Figure \ref{heat_N_numerical1_t} in the form of dissipation near the boundary. With $h = 0.01$, this drop in accuracy is illustrated in the error plot of Figure \ref{heat_N_UTM1_HL}.
	
	\begin{remark}
	We may consider different spacial discretizations of the heat equation in the IBVPs \eqref{heat_prob} or \eqref{heat_prob_N}. For example, the standard forward one-sided discretization of the heat equation,
	\begin{equation}
		\dot{q}_n(t) = \frac{q_{n+2}(t) - 2q_{n+1}(t) + q_n(t)}{h^2},
		\label{heat_forward}
	\end{equation}
	gives rise to the dispersion relation
	\begin{equation}
		W(k) = \frac{2 e^{ikh} - e^{2ikh} - 1}{h^2},
		\label{W_heat_forward}
	\end{equation}
	and  
	\begin{equation}
		q_n(T) = \frac{1}{2 \pi} \int_{-\pi/h}^{\pi/h} e^{iknh} e^{-WT}\hat{q}(k,0)\,dk + \frac{1}{2 \pi} \int_{-\pi/h}^{\pi/h} e^{ik(n+1)h} e^{-WT}\left[ \frac{\left(2 - e^{ikh}\right) f_{0} - f_{1}}{h}  \right]\,dk,
		\label{soln1_heat_forward}
	\end{equation} 
	using a forward discrete Fourier transform that starts at $n = 0$. Regardless of the starting index and available boundary conditions from the continuous problem, the second integral in the ``solution'' has zero contribution, \textit{i.e.}, the solution does not depend on any boundary information at all. This can be done by deforming off the real line as in Section \ref{advec_forward_halfline}, since the dispersion relation \eqref{W_heat_forward} with $z = e^{ikh}$ has all nonnegative degrees and $e^{-WT}$ is bounded in the upper half-plane. Thus, \eqref{heat_forward} gives rise to an ill-conditioned semi-discrete problem, relative to its continuous counterpart. 
	
	A similar issue arises when we consider a backward one-sided discretization for $q_{xx}$, except here the dispersion relation 
	\begin{equation}
		W(k) = \frac{2 e^{-ikh} - e^{-2ikh} - 1}{h^2},
		\label{W_heat_backward}
	\end{equation}
	does not permit the removal of either unknown from  
	\begin{equation}
		\hspace{-25pt}q_n(T) = \frac{1}{2 \pi} \int_{-\pi/h}^{\pi/h} e^{iknh} e^{-WT} \hat{q}(k,0) \,dk +  \frac{1}{2 \pi} \int_{-\pi/h}^{\pi/h} e^{iknh} e^{-WT} \left[ \frac{f_{-2} - \left(2 - e^{-ikh}\right) f_{-1}}{h} \right]\,dk.
		\label{soln_heat_backward}
	\end{equation}
	Although the discretization is first-order accurate, it has a second-order stencil with a dispersion relation, which has a nontrivial symmetry. Even so, it is not feasible to deform to the region where the global relation with this symmetry is valid, see Figure \ref{symm_heat_backward}.
	\begin{figure}[tb]
		\begin{center}
			\includegraphics[width=0.4\linewidth]{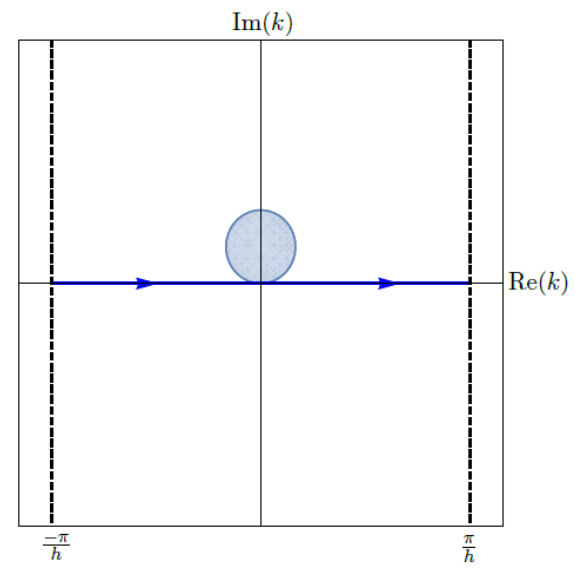}
			\caption{The shaded region depicts where the global relation of the backward one-sided discretization of $q_{xx}$ is valid, with $ k \rightarrow \nu_1(k)$ from \eqref{W_heat_backward}.}
			\label{symm_heat_backward}
		\end{center}
	\end{figure}
	Thus, this one-sided discretization is also problematic, requiring too much information from boundary nodes. 
	\end{remark}


	\subsection{Higher-Order Discretization of $\bms{q_t = q_{xx}}$ with Dirichlet boundary condition}
\begin{figure}[tb]
	\raggedleft
	\begin{subfigure}[t]{.45\textwidth}
			\centering
  			\includegraphics[width=0.9\linewidth]{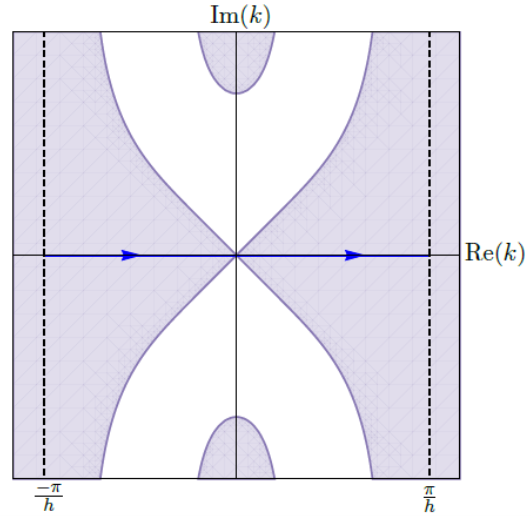}
  			\caption{}
			\label{W_heat_centered4_fig}
	\end{subfigure}\hfill 
	\begin{subfigure}[t]{.45\textwidth}
			\centering
  			\includegraphics[width=1.04\linewidth]{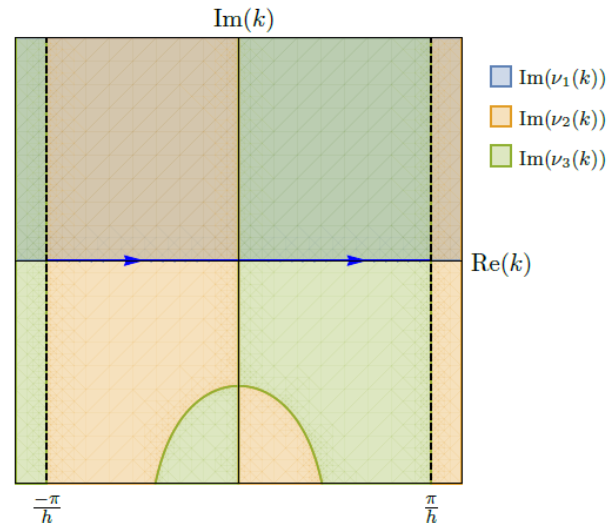}
  			\caption{}
			\label{symm_heat_centered4_fig}
	\end{subfigure}	
	\caption{(a) The shaded regions depict where $\text{Re}(-W) \leq 0$ and $e^{-WT}$ is bounded with the dispersion relation \eqref{W_heat_centered4}. (b) The shaded regions depict where the global relation with $ k \rightarrow \nu_i(k)$ is valid from \eqref{W_heat_centered4}.} 
	\label{heat_centered4_fig}
\end{figure}

	As in Section \ref{advec2_highorder_halfline}, we can apply higher-order discretizations to the heat equation where the nontrivial symmetries are not enough to eliminate unknowns. Consider the heat equation in \eqref{heat_prob} with the standard centered fourth-order discretization:
	\begin{equation}
		\dot{q}_n(t) =\frac{-q_{n-2}(t) + 16 q_{n-1}(t) - 30q_{n}(t) + 16q_{n+1}(t) - q_{n+2}(t)}{12h^2}.
		\label{heat_centered4}
	\end{equation}
	After several tedious steps, the global relation is
		\begin{align}
		e^{WT} \hat{q}(k,T) - \hat{q}(k,0) - F(k,T)&= 0, \quad \text{Im}(k) \leq 0,
		\label{GR_heat_centered4}
	\end{align}
	where
	$$F(k,T) = \frac{ - e^{-ikh} f_{-1} + 16 e^{-ikh} f_{0} - e^{-2ikh} f_{0} - 16 f_1 + e^{ikh} f_{1} + f_2}{12 h},$$
	with dispersion relation
	\begin{equation}
		W(k) = \frac{e^{-2ikh} - 16 e^{-ikh} + 30 - 16e^{ikh} + e^{2ikh}}{12h^2}.
		\label{W_heat_centered4}
	\end{equation}
	Solving for $\hat{q}(k,T)$ and taking the inverse transform, we obtain
	\begin{align}\begin{split}
		q_n(T) &= \frac{1}{2\pi} \int_{-\pi/h}^{\pi/h} e^{iknh} e^{-WT} \hat{q}(k,0)\,dk + \frac{1}{2\pi} \int_{-\pi/h}^{\pi/h} e^{iknh} e^{-WT} F(k,T)\,dk.
		\label{soln1_heat_centered4}
	\end{split}\end{align}
	Since we are given the Dirichlet boundary condition, $f_{-1}(W,T)$, $f_1(W,T)$, and $f_2(W,T)$ are unknown and must be removed from \eqref{soln1_heat_centered4}. The dispersion relation gives the nontrivial symmetries 
	\begin{align*}
		\nu_1(k) &= -k, \\
		\nu_2(k)& = \frac{i}{h} \ln \left(\frac{e^{-i k h}}{2}  \left[16 e^{i k h}-e^{2 i k h} - 1 + \sqrt{\left(-16 e^{i k h}+e^{2 i k h}+1\right)^2-4 e^{2 i k h}}\right]\right), \\
		\nu_3(k) &= \frac{i}{h} \ln \left(\frac{e^{-i k h}}{2}  \left[16 e^{i k h}-e^{2 i k h} - 1 - \sqrt{\left(-16 e^{i k h}+e^{2 i k h}+1\right)^2-4 e^{2 i k h}}\right]\right),
	\end{align*} 
	where the branch cut for the square-root function is chosen to be on the positive real line. Figure \ref{symm_heat_centered4_fig} depicts where in the complex $k$-plane the global relation \eqref{GR_heat_centered4} with $k \rightarrow \nu_i$ is valid, while Figure \ref{W_heat_centered4_fig} shows the shaded regions as decay due to $e^{-WT}$. Both figures imply that the integration path on the real line need not be deformed to use all three nontrivial symmetries. Even so, Figure \ref{symm_heat_centered4_fig} tells us that both $\nu_{2,3}(k)$ may only be used to remove one unknown, since there is no region where both symmetries are simultaneously valid. Hence, the three symmetries can only remove two unknowns and we must introduce a fourth equation.	
	With the given Dirichlet boundary condition, the heat equation itself gives a second-derivative boundary condition: 
	\begin{align*}
		q_{xx}(0,t) = q_{t}(0,t) = \frac{\partial}{\partial t} q(0,t) = u'(t) = v(t).
	\end{align*}
	To not introduce further unknowns, we must discretize $q_{xx}(0,t)$ using the standard centered stencil, so that
	\begin{align*}
		\frac{q_{-1}(t) - 2 q_{0}(t) + q_{1}(t)}{h^2} = v(t) \quad \Rightarrow \quad \frac{f_{-1} - 2 f_{0} + f_{1}}{h^2} = V, \quad\quad V(W,T) = \int_0^T e^{Wt} v(t)\, dt.
	\end{align*}
	This discretization drops the solution in accuracy from $\mathcal{O}(h^4)$ to $\mathcal{O}(h^2)$, but it allows a well-posed semi-discrete solution. The global relation with $ k \rightarrow \nu_i(k)$ and this discretized boundary condition give the system of equations
	\[\begin{dcases}
		0 = e^{WT} \hat{q}\left(-k,T\right) - \hat{q}\left(-k,0\right) - F\left(-k,T\right), \quad& \text{Im}(k) \geq 0, \\
		0 = e^{WT} \hat{q}\left(\nu_2,T\right) - \hat{q}\left(\nu_2,0\right) - F\left(\nu_2,T\right), & \text{Im}(\nu_2) \leq 0, \\
		0 = e^{WT} \hat{q}\left(\nu_3,T\right) - \hat{q}\left(\nu_3,0\right) - F\left(\nu_3,T\right), & \text{Im}(\nu_3) \leq 0, \\
		h^2 V = f_{-1} - 2 f_{0} + f_{1}, & k \in \mathbb{C}.
	\end{dcases}\]
	Solving for the unknowns, we find
	\begin{align}
		q_n(T) &= \frac{1}{2\pi} \int_{-\pi/h}^{\pi/h} e^{iknh} e^{-WT} \left[\hat{q}(k,0) - \hat{q}(-k,0)\right]\,dk + \frac{1}{2\pi} \int_{-\pi/h}^{\pi/h} e^{iknh} e^{-WT} \tilde{F}(k,T)\,dk,
		\label{soln_heat_centered4}
	\end{align}
	after deforming away the integral with $\hat{q}(-k,T)$, where
	$$\tilde{F}(k,T) = \frac{e^{-2 i k h} \left(14 e^{i k h}-14 e^{3 i k h}+e^{4 i k h}-1\right)}{12 h} f_0 + \frac{ h e^{-i k h} \left(e^{2 i k h}-1\right)}{12} V.$$
	Note that $\tilde{F}(k,T)$ has no dependence on $\nu_{2,3}$. As $h \rightarrow 0$, the semi-discrete solution \eqref{soln_heat_centered4} converges to \eqref{soln_heat_cont} and $\lim_{h\rightarrow 0} \tilde{F} = -2 i k F_0$. Thus, the semi-discrete solution correctly loses dependence on the second-derivative boundary condition in the continuum limit.


\section{The Linear Schr\"{o}dinger Equation}

	We consider the linear Schr\"{o}dinger (LS) equation
	\begin{equation}
		i q_t + \frac{1}{2} q_{xx} = 0 \quad\text{ or }\quad q_t = \frac{i}{2}q_{xx}.
		\label{LS}
	\end{equation}
	In contrast to the dissipative heat equation, this problem is dispersive.


\subsection{Centered Discretization of $\bms{q_t = \frac{i}{2}q_{xx}}$ with Dirichlet boundary condition}
	We begin with the half-line IBVP
	\begin{equation}\begin{dcases}
		q_t = \tfrac{i}{2} q_{xx},& x > 0,\, t > 0, \\
		q(x,0) = \phi(x),& x > 0,\\
		q(0,t) = u(t),& t > 0,
		\label{LS_prob}
	\end{dcases}\end{equation}
	using a centered discretization for $q_{xx}$,
	\begin{equation}		
		\dot{q}_n(t) = \frac{i}{2}\left( \frac{q_{n+1}(t) - 2 q_n(t) + q_{n-1}(t)}{h^2}\right).
		\label{LS_centered_D}
	\end{equation}
	The local and dispersion relations are, respectively, 
	\begin{align}
		\partial_t \left(e^{-iknh} e^{Wt} q_n \right) &= \frac{i}{2 h^2}\Delta \left(e^{-ik(n-1)h} e^{Wt} q_{n} - e^{-iknh} e^{Wt} q_{n-1} \right),
		\label{LR_LS_centered}
	\end{align} 
	\begin{equation}
		W(k) = \frac{i}{2}\left( \frac{2 - e^{ikh} - e^{-ikh}}{h^2}\right) = \frac{i}{h^2} \left[ 1 - \cos(kh) \right].
		\label{W_LS_centered}
	\end{equation}
	With the Dirichlet boundary condition, our transforms begin at $n = 1$ instead of at $n = 0$, resulting in the global relation
	\begin{align}
		e^{WT} \hat{q}(k,T) - \hat{q}(k,0) - \frac{i}{2}\left[ \frac{ e^{-ikh} f_{0} - f_1}{h} \right] &= 0, \quad \text{Im}(k) \leq 0.
		\label{GR_LS_centered}
	\end{align}
	To obtain our ``solution'' formula, we take the inverse transform,
	\begin{align}
		q_n(T) &= \frac{1}{2 \pi} \int_{-\pi/h}^{\pi/h} e^{iknh} e^{-WT}\hat{q}(k,0)\,dk + \frac{1}{2 \pi} \int_{-\pi/h}^{\pi/h}\frac{ i e^{iknh} e^{-WT}}{2}\left[ \frac{ e^{-ikh} f_{0} - f_1}{h} \right]\,dk.
		\label{soln1_LS_centered}
	\end{align}
	
	The dispersion relation \eqref{W_LS_centered} admits the nontrivial symmetry $\nu_1 (k) = -k$ up to periodic copies. Hence, the global relation \eqref{GR_LS_centered} with $\nu_1(k)$ is valid in the upper-half plane, including the real line, so that there is no need to deform in order to eliminate the unknown $f_{1}(W,T)$ in \eqref{soln1_LS_centered}. We find
	$$f_1 = e^{i k h} f_0 - 2 i h \left[\hat{q}(-k,0) - e^{T W} \hat{q}(-k,T)\right],$$
	so that
	\begin{align}
	\begin{split}
		q_n(T) &= \frac{1}{2 \pi} \int_{-\pi/h}^{\pi/h} e^{iknh} e^{-WT}\hat{q}(k,0)\,dk + \frac{1}{2 \pi} \int_{-\pi/h}^{\pi/h} e^{iknh} e^{-WT} \left[\frac{i \left(e^{-i k h} - e^{i k h}\right)}{2 h} f_{0} - \hat{q}(-k,0) \right]\,dk \\
		&\quad\, + \frac{1}{2 \pi} \int_{-\pi/h}^{\pi/h} e^{iknh}  \hat{q}(-k,T)\,dk.
		\label{soln2_LS_centered}
	\end{split}
	\end{align}
	As before, one shows that the last term does not contribute. Therefore, the solution to this half-line IBVP is
	 \begin{align}
	\begin{split}
		q_n(T) &= \frac{1}{2 \pi} \int_{-\pi/h}^{\pi/h} e^{iknh} e^{-WT}\hat{q}(k,0)\,dk  - \frac{1}{2 \pi} \int_{-\pi/h}^{\pi/h} e^{iknh} e^{-WT} \left[ \hat{q}(-k,0) - \frac{\sin(kh)}{h} f_{0} \right]\,dk.
		\label{soln_LS_centered}
	\end{split}
	\end{align}
	
	Since a centered stencil was used to obtain \eqref{soln_LS_centered}, the modified PDE which this semi-discrete solution better approximates solutions of is similar to \eqref{heat_centered_modified_eqn} derived for the heat equation. Instead of being dissipative, we have the dispersive PDE $p_t = (i/2)p_{xx} + (i h^2/24) p_{4x}$. Solution \eqref{soln_LS_centered} solves this modified PDE to fourth-order, where the dispersive behavior of the second term on the right-hand side is evident for large $t$ and fixed $h$. Because this term is $\mathcal{O}\left(h^2\right)$, we can diminish the excess dispersion by decreasing $h$. 
	
	The continuous UTM solution \cite{bernard_fokas} to \eqref{LS_prob} is	
	\begin{align}
		q(x,T) &= \frac{1}{2 \pi} \int_{-\infty}^{\infty} e^{ikx} e^{-\tilde{W}T} \hat{q}(k,0)\, dk - \frac{1}{2 \pi} \int_{-\infty}^{\infty} e^{ikx} e^{-\tilde{W}T} \left[\hat{q}(-k,0) - k F_0 \right]\,dk,
		\label{soln_LS_cont}
	\end{align}
	with dispersion relation $\tilde{W}(k) = i k^2/2$ and nontrivial symmetry $\tilde{\nu}_1 = -k$. The semi-discrete solution \eqref{soln_LS_centered} converges to its continuous counterpart solution \eqref{soln_LS_cont} in the continuum limit.

\begin{figure}[h!]	
		\raggedleft
		\begin{subfigure}[t]{.45\textwidth}
			\centering
  			\includegraphics[width=1\linewidth]{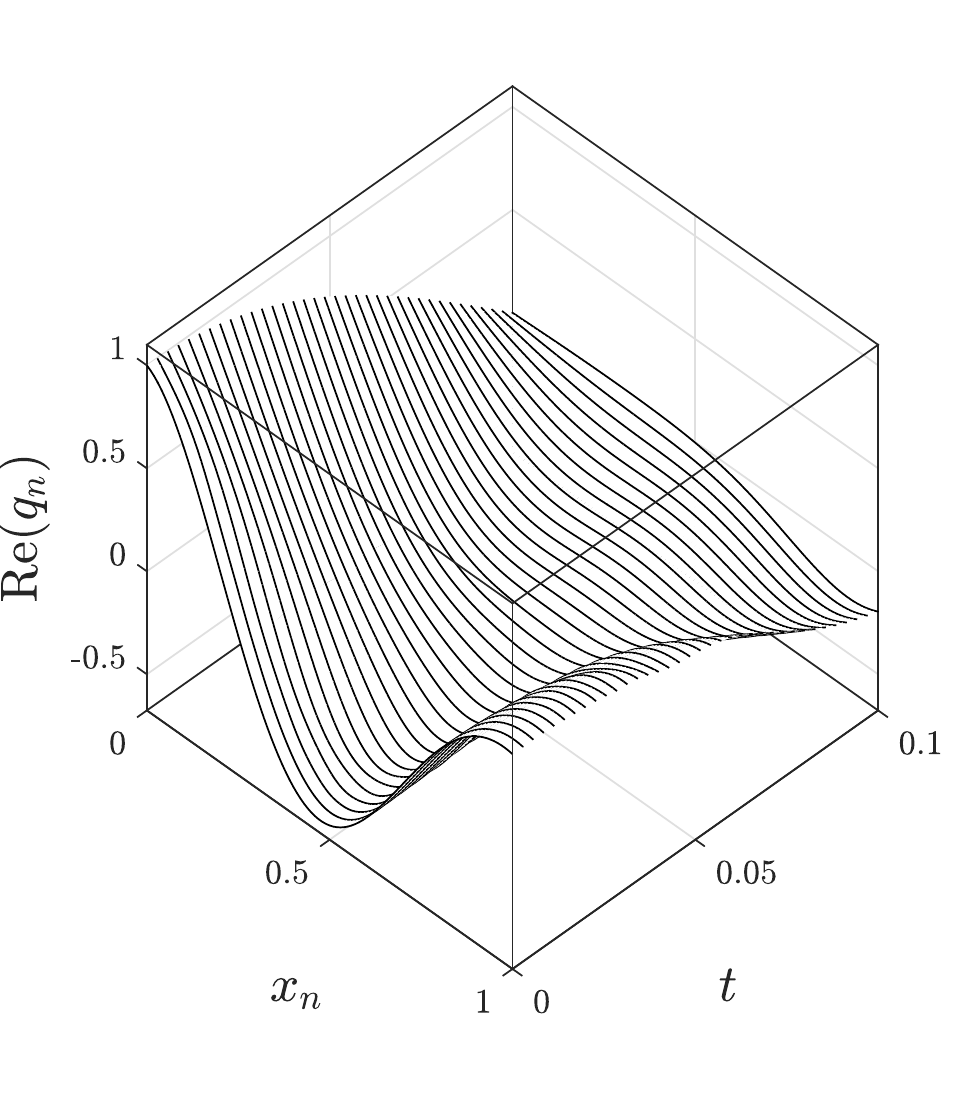}
  			\caption{}
  			\label{}
		\end{subfigure}\hfill 
		\begin{subfigure}[t]{.45\textwidth}
			\centering
  			\includegraphics[width=1\linewidth]{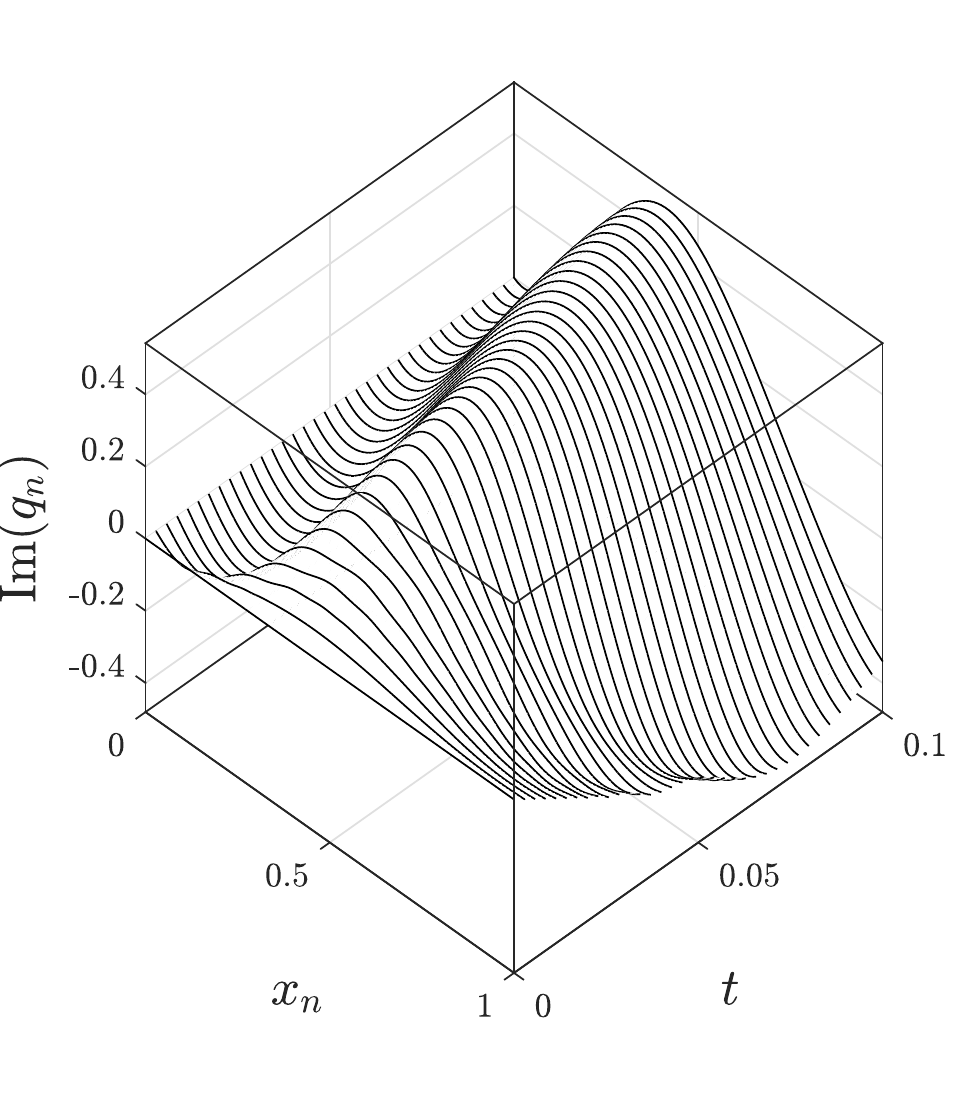}
  			\caption{}
  			\label{}
		\end{subfigure}
		\begin{subfigure}[t]{.45\textwidth}
			\centering
  			\includegraphics[width=1\linewidth]{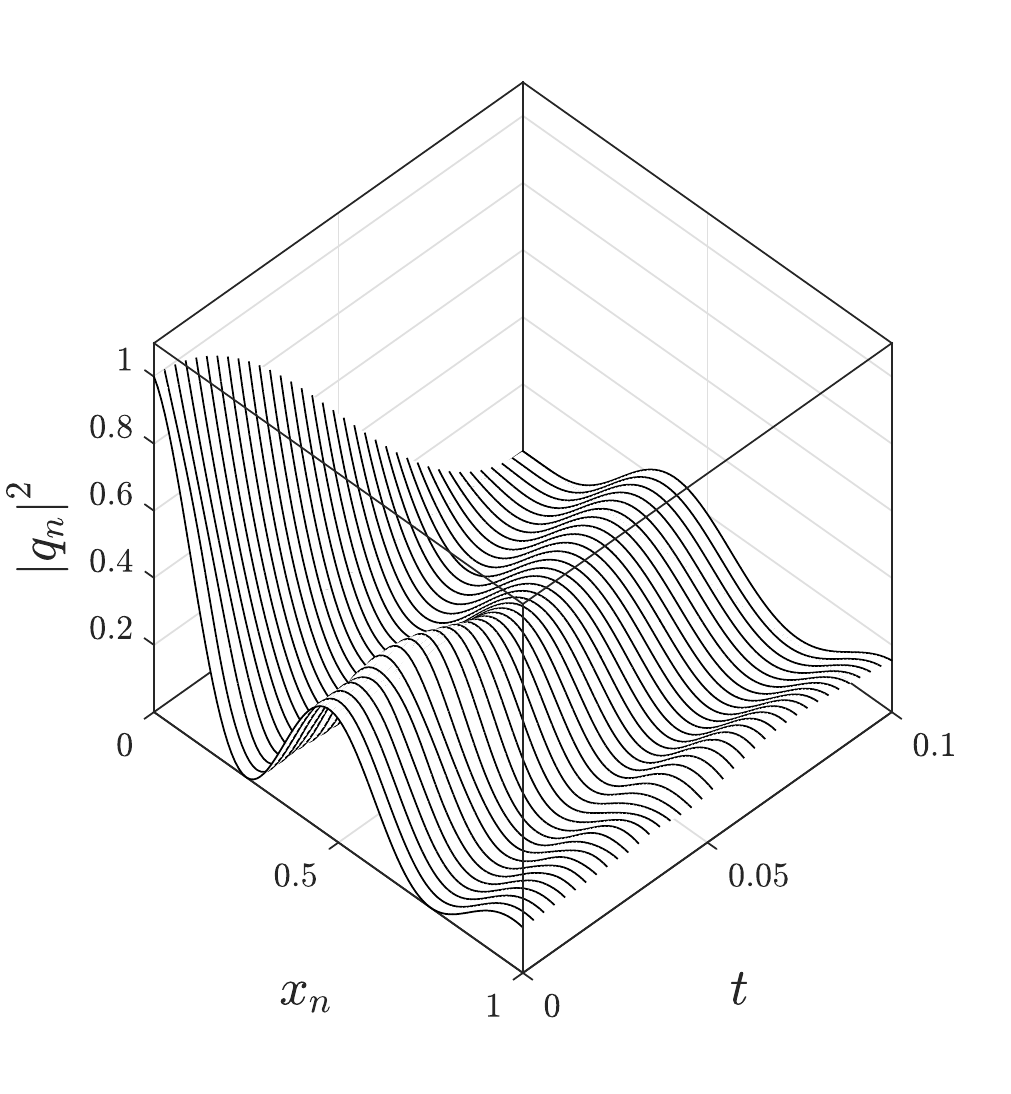}
  			\caption{}
  			\label{}
		\end{subfigure}\hfill
		\begin{subfigure}[t]{.45\textwidth}
			\centering
  			\includegraphics[width=1.05\linewidth]{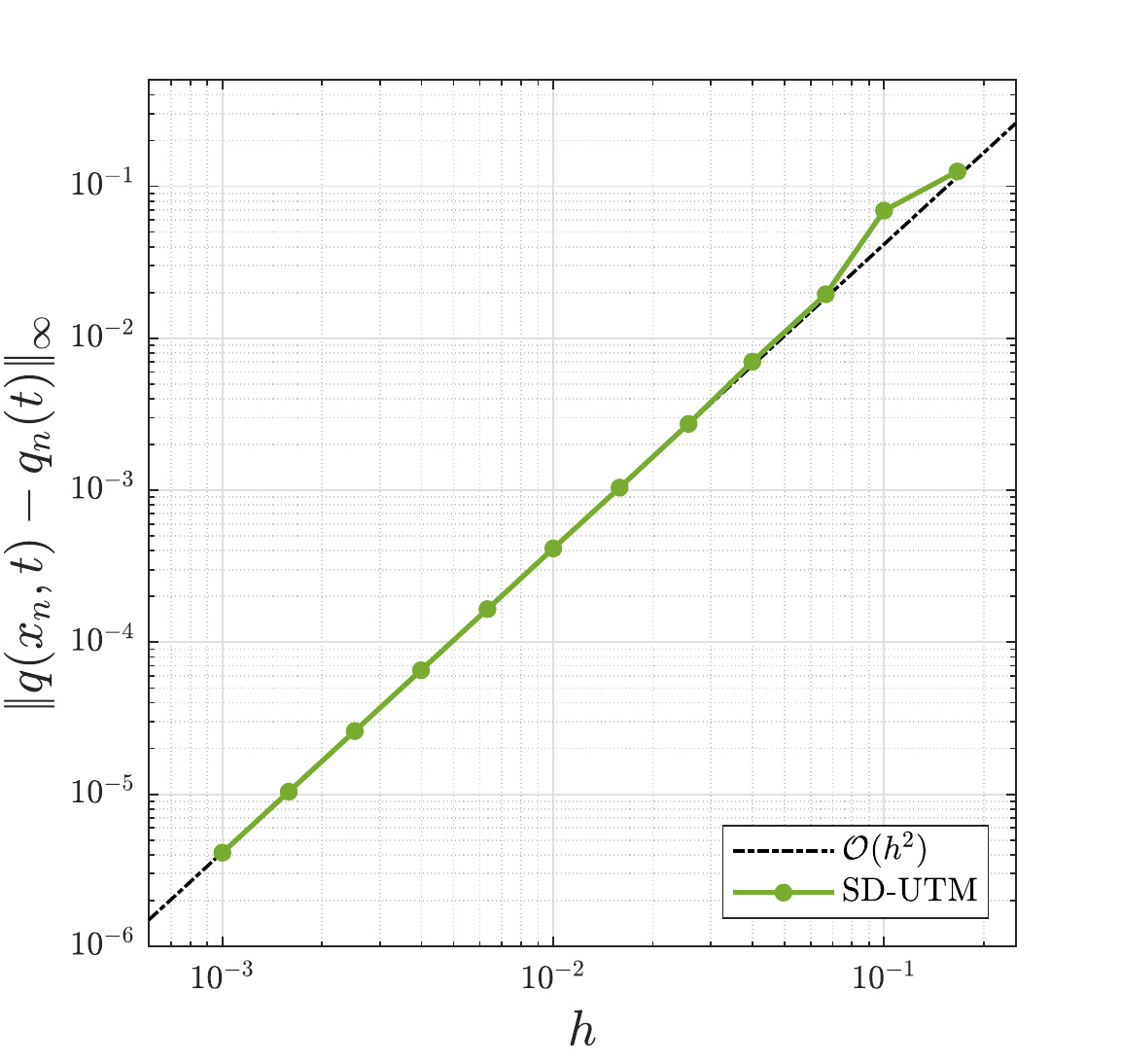}
  			\caption{}
  			\label{}
		\end{subfigure}	
		\caption{(a) - (c) Real and imaginary parts and modulus squared of the semi-discrete solution profiles \eqref{soln_LS_centered} at various $t$ for IBVP \eqref{LS_numerical1_HL} with $h = 0.01$. (d) Error plot of the semi-discrete solution \eqref{soln_LS_centered} relative to the exact solution as $h \rightarrow 0$ with $t = 0.1$.}
		\label{LS_UTM1_HL}
	\end{figure}

	We examine the numerical solution to
\begin{equation}
\begin{dcases}
	q_t = \tfrac{i}{2} q_{xx},& x > 0,\, t > 0, \\
	q(x,0) = \phi(x) = e^{-x} \cos \left(2 \pi x \right),& x > 0,  \\
	q(0,t) = u(t) = \cos \left(5 \pi t \right),& t > 0,
\end{dcases}
\label{LS_numerical1_HL}
\end{equation}	
	Like the heat equation, the continuous solution to this problem can be written in terms of error functions of imaginary argument. The semi-discrete solution for the second-order finite-difference approximation \eqref{LS_centered_D} is given by \eqref{soln_LS_centered}. Figure \ref{LS_UTM1_HL} shows the dispersive nature of the real and imaginary components of the solution, along with the square of the modulus. 


\subsection{Centered Discretization of $\bms{q_t = \frac{i}{2}q_{xx}}$ with Neumann boundary condition} \label{LS_neumann_halfline}
	We consider the same centered-discretized LS equation as above with a Neumann boundary condition:
	\begin{equation}\begin{dcases}
		q_t = \tfrac{i}{2} q_{xx},& x > 0,\, t > 0, \\
		q(x,0) = \phi(x),& x > 0,\\
		q_x(0,t) = v(t),& t > 0.
		\label{LS_prob_N}
	\end{dcases}\end{equation}
	Here, $q_0(t)$ is unknown, and we choose the discrete Fourier transform to start from $n= 0$ instead of from $n = 1$. The local and dispersion relations, \eqref{LR_LS_centered} and \eqref{W_LS_centered} respectively, remain unchanged. The global relation is 
	\begin{align}
		e^{WT} \hat{q}(k,T) - \hat{q}(k,0) - \frac{i}{2}\left[ \frac{f_{-1} - e^{ikh} f_{0} }{h} \right] &= 0, \quad \text{Im}(k) \leq 0.
		\label{GR_LS_centered_N}
	\end{align}
	Using the inverse transform,
	\begin{align}
		q_n(T) &= \frac{1}{2 \pi} \int_{-\pi/h}^{\pi/h} e^{iknh} e^{-WT}\hat{q}(k,0)\,dk + \frac{1}{2 \pi} \int_{-\pi/h}^{\pi/h}\frac{ i e^{iknh} e^{-WT}}{2}\left[ \frac{f_{-1} - e^{ikh} f_{0} }{h} \right]\,dk.
		\label{soln1_LS_centered_N}
	\end{align}
	Like the heat equation with a Neumann boundary condition at $x = 0$, we apply the standard backward discretization to $q_x(0,t)$ so as not to introduce new unknowns,
	$$\frac{q_0(t) - q_{-1}(t)}{h} = v(t).$$
	As discussed in Section \ref{neumann_halfline}, this discretization drops the accuracy to $\mathcal{O}(h)$, visible in solution profiles as dissipation near the boundary. The global relation \eqref{GR_LS_centered_N} and the time transform of the discretized boundary condition give the system
	\[\begin{dcases*}
		e^{WT} \hat{q}(-k,T) - \hat{q}(-k,0) - \frac{i}{2}\left[ \frac{f_{-1} - e^{-ikh} f_{0} }{h} \right] = 0, \\
		\frac{f_0 - f_{-1}}{h} = V(t),
	\end{dcases*}\]
	for the two unknowns $f_{-1}(W,T)$ and $f_0(W,T)$, where
	$$V(W,T) = \int_0^T e^{Wt} v(t)\, dt.$$
	Solving the system and substituting into \eqref{soln1_LS_centered_N}, we have
	\begin{align}
	\begin{split}
		\hspace{-15pt}q_n(T) &= \frac{1}{2 \pi} \int_{-\pi/h}^{\pi/h} e^{iknh} e^{-WT}\hat{q}(k,0)\,dk - \frac{1}{2 \pi} \int_{-\pi/h}^{\pi/h} e^{iknh} e^{-WT} \left[ \frac{i\left( e^{i k h} + 1 \right)}{2} V(t) - e^{i k h} \hat{q}(-k,0)  \right]\,dk,
		\label{soln_LS_centered_N}
	\end{split}
	\end{align}
	after applying similar techniques as before to remove the integral term depending on $\hat{q}(-k,T)$. This limits to the solution \cite{bernard_fokas} of the continuous problem
	\begin{align}
		q(x,T) &= \frac{1}{2 \pi} \int_{-\infty}^{\infty} e^{ikx} e^{-\tilde{W}T} \hat{q}(k,0)\, dk - \frac{1}{2 \pi} \int_{-\infty}^{\infty} e^{ikx} e^{-\tilde{W}T} \left[i F_1 - \hat{q}(-k,0) \right]\,dk.
		\label{soln_LS_cont_N}
	\end{align}

\begin{figure}[h!]	
		\raggedleft
		\begin{subfigure}[t]{.45\textwidth}
			\centering
  			\includegraphics[width=1\linewidth]{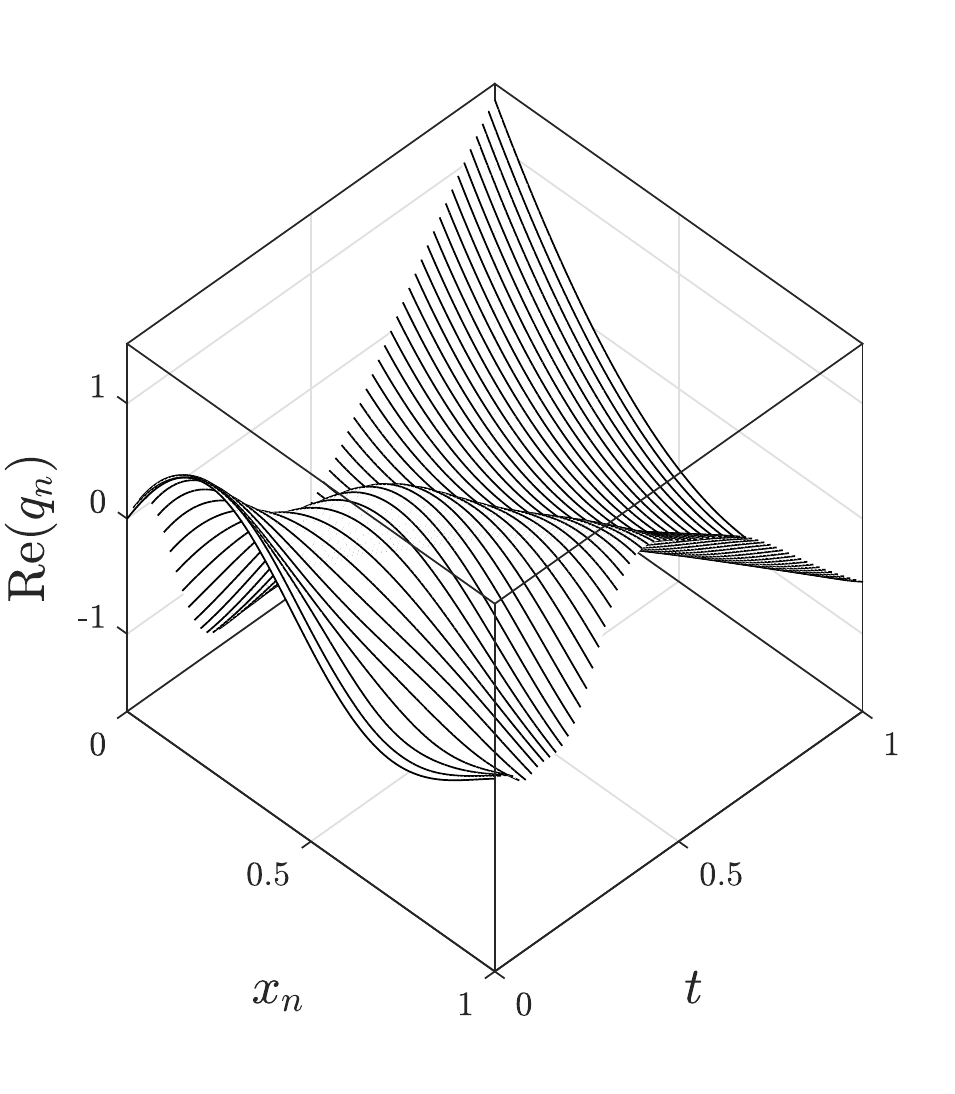}
  			\caption{}
  			\label{}
		\end{subfigure}\hfill 
		\begin{subfigure}[t]{.45\textwidth}
			\centering
  			\includegraphics[width=1\linewidth]{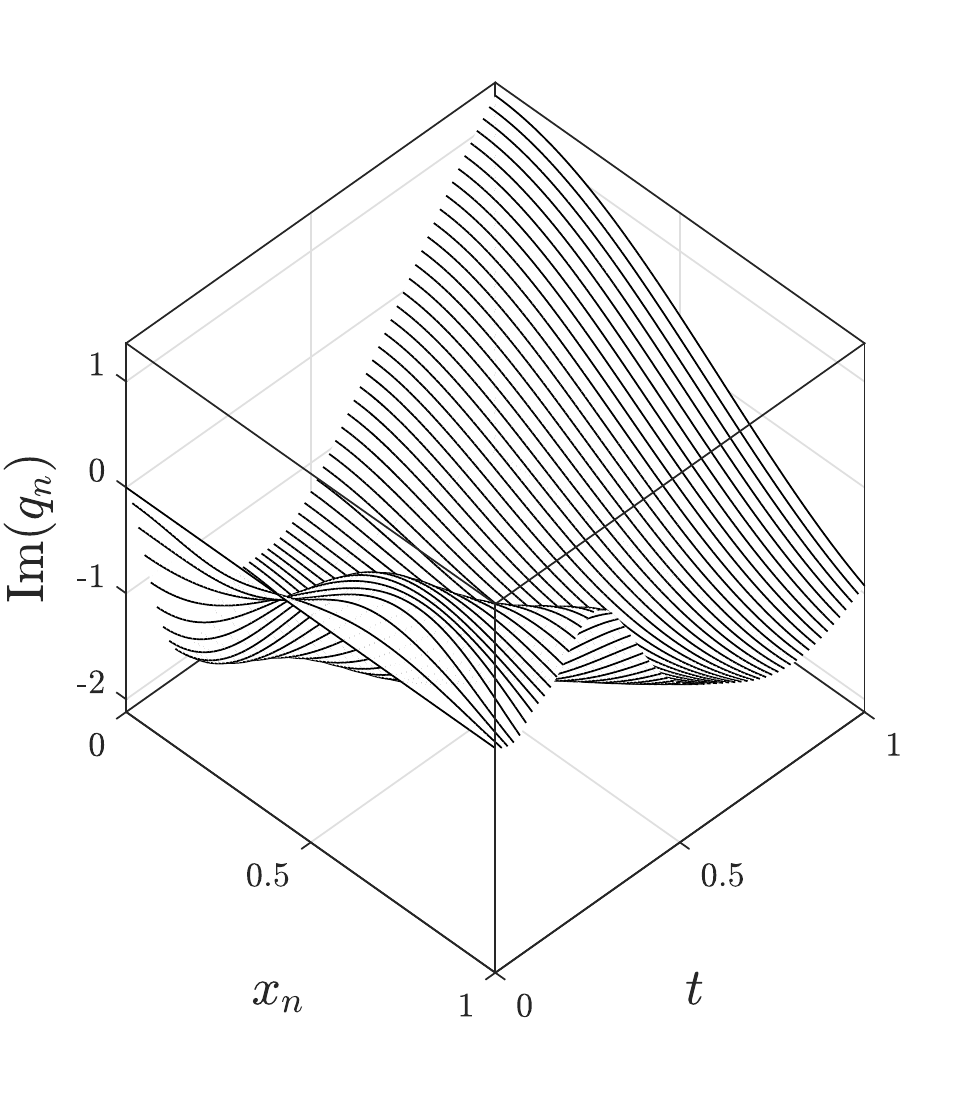}
  			\caption{}
  			\label{}
		\end{subfigure}
		\begin{subfigure}[t]{.45\textwidth}
			\centering
  			\includegraphics[width=1\linewidth]{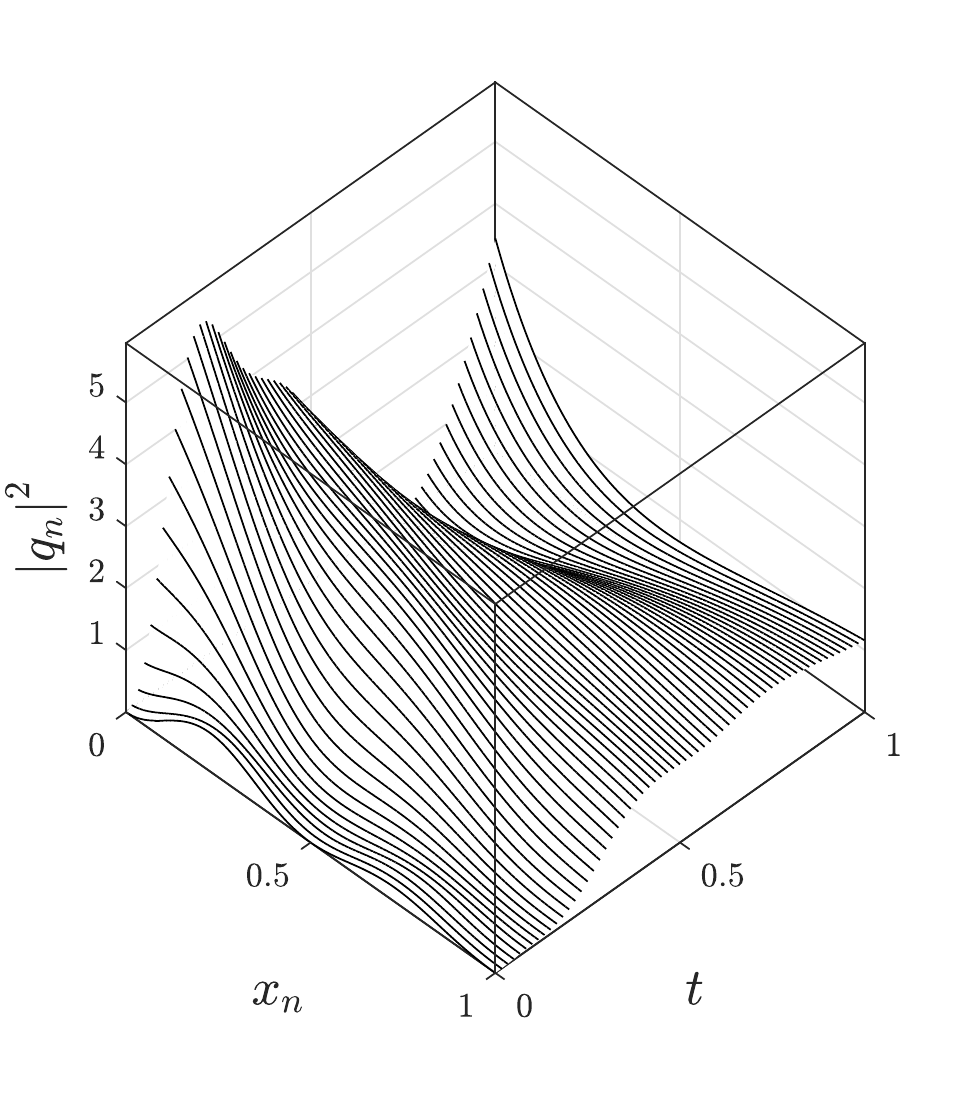}
  			\caption{}
  			\label{}
		\end{subfigure}\hfill
		\begin{subfigure}[t]{.45\textwidth}
			\centering
  			\includegraphics[width=1.05\linewidth]{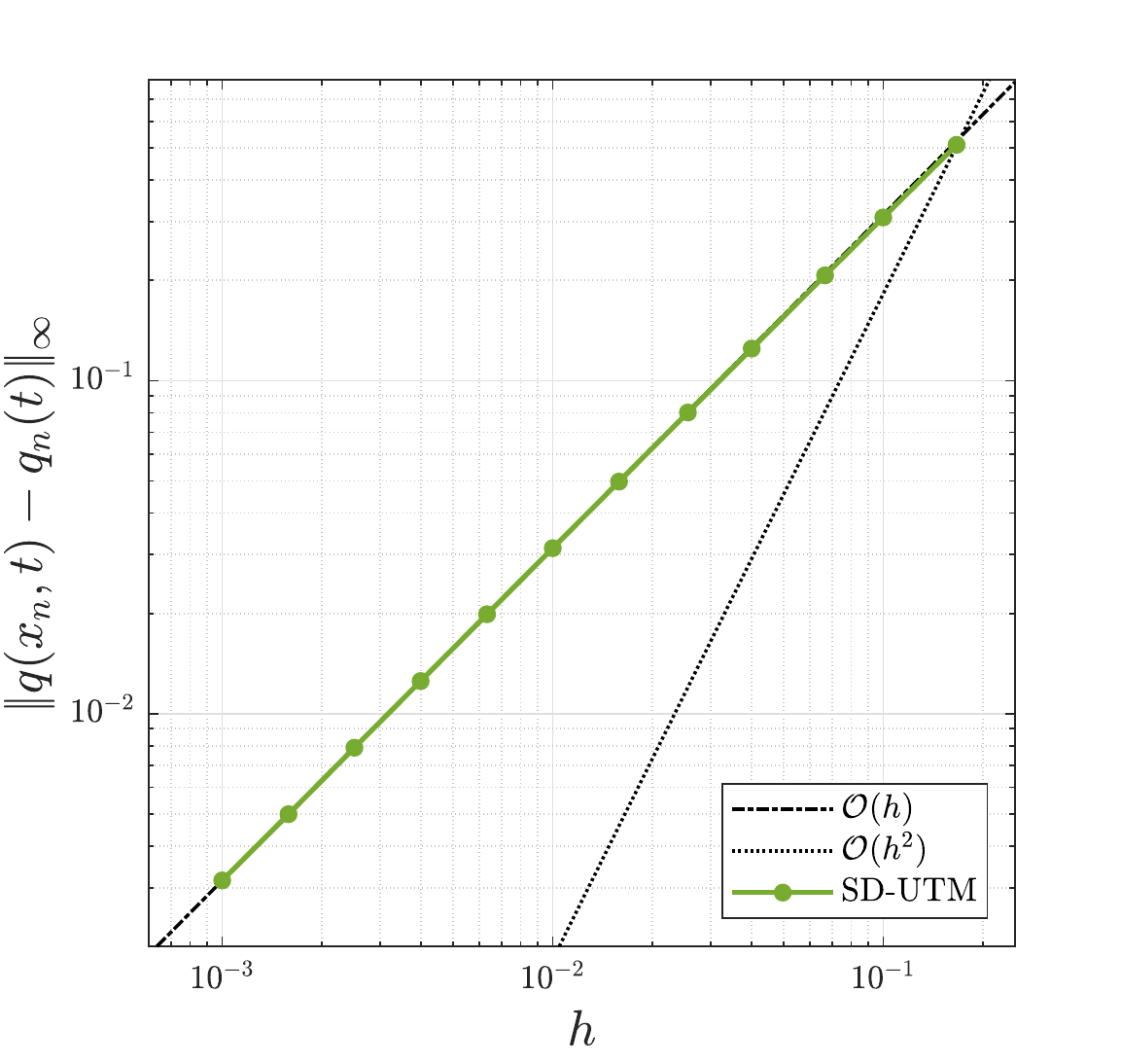}
  			\caption{}
  			\label{}
		\end{subfigure}	
		\caption{(a) - (c) Real and imaginary parts and modulus squared of the semi-discrete solution profiles \eqref{soln_LS_centered_N} at various $t$ for IBVP \eqref{LS_N_numerical1_HL} with $h = 0.01$. (d) Error plot of the semi-discrete solution \eqref{soln_LS_centered_N} relative to the exact solution as $h \rightarrow 0$ with $t = 1$.}
		\label{LS_N_UTM1_HL}
	\end{figure}

	Lastly, consider
\begin{equation}
\begin{dcases}
	q_t = \tfrac{i}{2} q_{xx},& x > 0,\, t > 0, \\
	q(x,0) = \phi(x) =  e^{-x} \sin \left(2 \pi x \right),& x > 0,  \\
	q_x(0,t) = v(t) = 2 \pi \cos \left(\pi t \right),& t > 0,
\end{dcases}
\label{LS_N_numerical1_HL}
\end{equation}	
with the Neumann condition discretized using the standard backward stencil, giving rise to the $\mathcal{O}(h)$ accurate solution \eqref{soln_LS_centered_N}, displayed in Figure \ref{LS_N_UTM1_HL}. 


\section{Small-Time Increments}\label{small_time_sec}

	For nonlinear IBVPs, the semi-discrete nor continuous UTM is applicable in general. Although the UTM can be used to solve IBVPs for integrable PDEs, our goal is broader: we are interested in numerically solving IBVPs for quasilinear PDEs \eqref{ibvp_eq}, where the most nonlocal stencil is applied to the linear problem. To do so accurately, we can employ split-step methods following the ideas from operator splitting. We rewrite the evolution PDE \eqref{ibvp_eq} as
	\begin{equation} 
		q_t = \mathcal{L}\left(q\right) + \mathcal{N}\left(q\right),
		\label{ibvp_eq_splitstep}
	\end{equation} 
	where $\mathcal{L}$ is a constant-coefficient linear differential operator of order $N$ and $\mathcal{N}$ is a nonlinear operator, both operators involving spatial derivatives of $q(x,t)$. The idea behind split-step methods is to separately solve the $N^{\text{th}}$-order linear IBVP, with
	\begin{equation}
		q_t = \mathcal{L}\left(q\right) = c \, q_{Nx},
		\label{eq_splitstep1}
	\end{equation} 
	and the nonlinear IBVP, with
	\begin{equation}
		q_t = \mathcal{N}\left(q\right) = F\left(q,q_x,\ldots,q_{(N-1)x} \right),
		\label{eq_splitstep2}
	\end{equation} 
	and combine them iteratively \cite{macnamara,randy_splitstep1,randy_splitstep2}. 
	
	The implementation of boundary conditions, other than periodic, is problematic for many numerical methods, including split-step methods. As with other finite-difference approaches, the use of high-degree spacial discretizations of the operators $\mathcal{L}$ and $\mathcal{N}$ might introduce ghost points or artificial boundary conditions that affect the overall performance of the numerical method. 
	
	We aim to overcome this problem by applying the semi-discrete UTM to the linear operator $\mathcal{L}$ in a split-step approach, while correctly incorporating boundary conditions. This split-step method requires the repeated computation of the solution to this linear problem with $t \ll 1$, but the integral representations from semi-discrete UTM can be expensive to compute. In what follows, we evaluate the semi-discrete UTM solutions using a $t \ll 1$ approximation to derive an approximate semi-discrete UTM solution with predetermined accuracy in $t$. Here, we demonstrate this procedure for the advection equation $q_t = -c\, q_x$ on the half line with the standard backward stencil \eqref{advec2_backward} applied to $q_x$, while more details and further investigations will be presented in a future paper. Since a split-step approach solves an updated IBVP starting from $t_0$, we generalize the original IBVP \eqref{advec2_prob} to
	\begin{equation}\begin{dcases}
		q_t = -c\,q_{x},& x > 0,\, t > t_0 \\
		q\left(x,t_0\right) = \phi(x),& x > 0,\\
		q(0,t) = u(t),& t > t_0,
		\label{advec2_prob_splitstep}
	\end{dcases}\end{equation}
	where $\phi(x)$ is the output from the previous step.
	
	Starting from $t_0$, the time transforms from the semi-discrete UTM are redefined as
	$$f_j\left(W,t_0,T\right) = \int_{t_0}^T e^{Wt} q_j(t) \, dt, \quad k \in \mathbb{C},$$
	which, for the IBVP	\eqref{advec2_prob_splitstep}, gives the global relation
	\begin{align}
		&&\hspace{-55pt} \sum_{n=1}^{\infty} h \int_{t_0}^T \left[ \partial_t \left(e^{-iknh} e^{Wt} q_n \right) + \frac{c}{h}\Delta \left(e^{-iknh} e^{Wt} q_{n-1} \right) \right] dt &= 0 \notag \\
		\Rightarrow &&\hspace{-55pt} e^{WT} \hat{q}(k,T) - e^{W t_0}\hat{q}\left(k,t_0\right) - c e^{-ikh} f_{0} &= 0,
		\label{GR_advec2_backward_transform_splitstep}
	\end{align}
	valid for $\text{Im}(k) \leq 0$. Solving for $\hat{q}(k,T)$ and inverting, we obtain
	\begin{align}\begin{split}
		q_n\left(T; t_0\right) &= \frac{1}{2\pi} \int_{-\pi/h}^{\pi/h} e^{iknh} e^{-W\left(T-t_0\right)} \hat{q}\left(k,t_0\right)\,dk + \frac{c}{2\pi} \int_{-\pi/h}^{\pi/h} e^{ik(n-1)h} e^{-WT} f_{0}\,dk.
		\label{soln_advec2_backward_splitstep}
	\end{split}\end{align}
	Following similar arguments as before, \eqref{soln_advec2_backward_splitstep} is the solution to the backward-discretized IBVP \eqref{advec2_prob_splitstep} with a given Dirichlet boundary condition at $x = 0$ and an initial condition at $t = t_0$. In what follows, we expand \eqref{soln_advec2_backward_splitstep} in $\tau = T - t_0 \ll 1$, around $\tau=0$, so as to obtain a convenient approximation to be used in a split-step method. 	
	
	We expand $e^{-W\tau}$ using its Taylor series about $\tau = 0$ up to arbitrary $\mathcal{O}\left(\tau^r\right)$, so that the integrals have polynomial dependence on time. For the first integral,
	\begin{align}
		\hspace{-25pt}\frac{1}{2\pi} \int_{-\pi/h}^{\pi/h} e^{iknh} e^{-W\tau} \hat{q}\left(k,t_0\right)\,dk &= \frac{1}{2\pi} \int_{-\pi/h}^{\pi/h} e^{iknh} \left[ 1 - W \tau + \frac{W^2}{2}\tau ^2  - \frac{ W^3}{6}\tau ^3 +\mathcal{O}\left(\tau ^4\right) \right] \hat{q}\left(k,t_0\right)\,dk \notag \\[11pt]
		\hspace{-25pt}&= \frac{1}{2\pi} \int_{-\pi/h}^{\pi/h} e^{iknh} \hat{q}\left(k,t_0\right)\,dk \,\,-\,\, \frac{\tau}{2\pi} \int_{-\pi/h}^{\pi/h} W e^{iknh}  \hat{q}\left(k,t_0\right)\,dk \notag\\
		\hspace{-25pt}&\quad\, + \frac{\tau^2}{4\pi} \int_{-\pi/h}^{\pi/h}  W^2 e^{iknh} \hat{q}\left(k,t_0\right)\,dk \,\,-\,\, \frac{\tau^3}{12\pi} \int_{-\pi/h}^{\pi/h} W^3 e^{iknh}  \hat{q}\left(k,t_0\right)\,dk \,\,+\,\, \mathcal{O}\left(\tau ^4\right)\notag\\[11pt]
		\begin{split}
		\hspace{-25pt}&= q_n\left(t_0\right) \,\,-\,\, \frac{\tau}{2\pi} \int_{-\pi/h}^{\pi/h} W e^{iknh}  \hat{q}\left(k,t_0\right)\,dk  + \frac{\tau^2}{4\pi} \int_{-\pi/h}^{\pi/h}  W^2 e^{iknh} \hat{q}\left(k,t_0\right)\,dk \notag\\
		\hspace{-25pt}&\quad\,- \,\, \frac{\tau^3}{12\pi} \int_{-\pi/h}^{\pi/h} W^3 e^{iknh}  \hat{q}\left(k,t_0\right)\,dk \,\,+\,\, \mathcal{O}\left(\tau ^4\right).\notag
	\end{split}	
	\end{align}
	
	The second integral of \eqref{soln_advec2_backward_splitstep} has time dependence in $e^{-WT}$ and $f_{0}\left(W,t_0,T\right)$. We consider these together:
	\begin{align*}
		e^{-WT} f_{0} &= e^{-WT} \int_{t_0}^T e^{Wt} q_0(t) \, dt = e^{-W\tau} \int_{0}^{\tau} e^{W\tilde{t}} q_0\left(\tilde{t} + t_0\right) \, d\tilde{t}.
	\end{align*}
	Since the limits of integration approach zero as $\tau \rightarrow 0$, we expand $ e^{W\tilde{t}} q_0\left(\tilde{t} + t_0\right)$ about $\tilde{t} = 0$. Up to third-order terms, 
	$$e^{-WT} f_{0} = q_0\left( t_0 \right) \tau + \frac{q_0'\left( t_0 \right) - W q_0\left( t_0 \right)}{2} \, \tau ^2  + \frac{q_0''\left( t_0 \right) - W q_0'\left( t_0 \right)+W^2 q_0\left( t_0 \right)}{6} \, \tau ^3 + \mathcal{O}\left(\tau^4\right)$$
	so that the second integral of \eqref{soln_advec2_backward_splitstep} reduces to
	\begin{align*}
		\frac{c}{2\pi} \int_{-\pi/h}^{\pi/h} e^{ik(n-1)h} e^{-WT} f_{0}\,dk &= \frac{c\, q_0\left( t_0 \right) \tau}{h} \delta_{1n} + \frac{c \tau^2}{4\pi} \int_{-\pi/h}^{\pi/h} e^{ik(n-1)h} \left[q_0'\left( t_0 \right)-W q_0\left( t_0 \right) \right]\,dk \\
		&\quad\, + \frac{c\tau^3}{12\pi} \int_{-\pi/h}^{\pi/h} e^{ik(n-1)h} \left[ q_0''\left( t_0 \right)- W q_0'\left( t_0 \right)+W^2 q_0\left( t_0 \right) \right]\,dk + \mathcal{O}\left(\tau^4 \right),
	\end{align*}
	where $\delta_{ij}$ is the Kronecker delta. Therefore, after defining $\phi(x) = q\left(x,t_0\right)$ and $u(t) = q(0,t)$, the solution \eqref{soln_advec2_backward_splitstep} up to third order in $\tau$ is
	\begin{align}
		q_n(\tau) &= \phi_n \,\,+\,\, K_1 \tau \,\, +\,\, K_2 \tau^2\,\, +\,\,K_3 \tau^3 \,\, +\,\, \mathcal{O}\left(\tau^4\right),
	\label{soln_advec2_splitstep}
	\end{align}
	with
	\begin{align*}
		K_1 (n) &= \frac{c\, u\left( t_0 \right)}{h} \delta_{1n} - \frac{1}{2\pi} \int_{-\pi/h}^{\pi/h} W e^{iknh}  \hat{\phi}\left(k\right)\,dk, \\
		K_2 (n) &= \frac{c }{4\pi} \int_{-\pi/h}^{\pi/h} e^{ik(n-1)h} \left[u'\left( t_0 \right)-W u\left( t_0 \right) \right]\,dk + \frac{1}{4\pi} \int_{-\pi/h}^{\pi/h}  W^2 e^{iknh} \hat{\phi}\left(k\right)\,dk, \\
		K_3 (n) &= \frac{c}{12\pi} \int_{-\pi/h}^{\pi/h} e^{ik(n-1)h} \left[ u''\left( t_0 \right)- W u'\left( t_0 \right)+W^2 u\left( t_0 \right) \right]\,dk -\frac{1}{12\pi} \int_{-\pi/h}^{\pi/h} W^3 e^{iknh}  \hat{\phi}\left(k\right)\,dk.
	\end{align*}
	A similar process can be repeated for other IBVPs.
	
\begin{figure}[tb]
	\raggedleft
	\begin{subfigure}[t]{.45\textwidth}
		\centering
  		\includegraphics[width=1\linewidth]{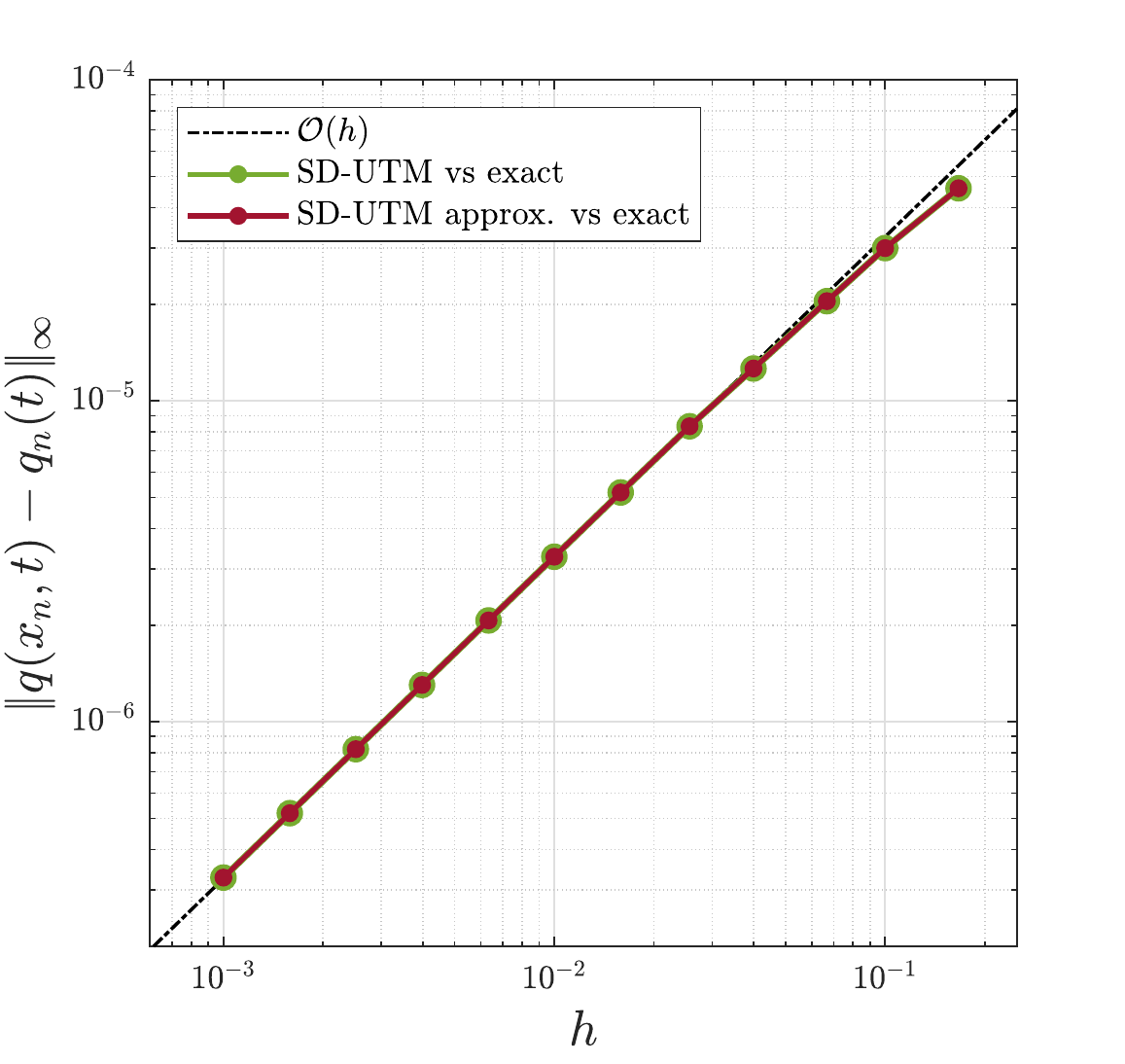}
  		\caption{}
  		\label{advec2_splitstep_a}
	\end{subfigure}\hfill
	\begin{subfigure}[t]{.45\textwidth}
		\centering
  		\includegraphics[width=1\linewidth]{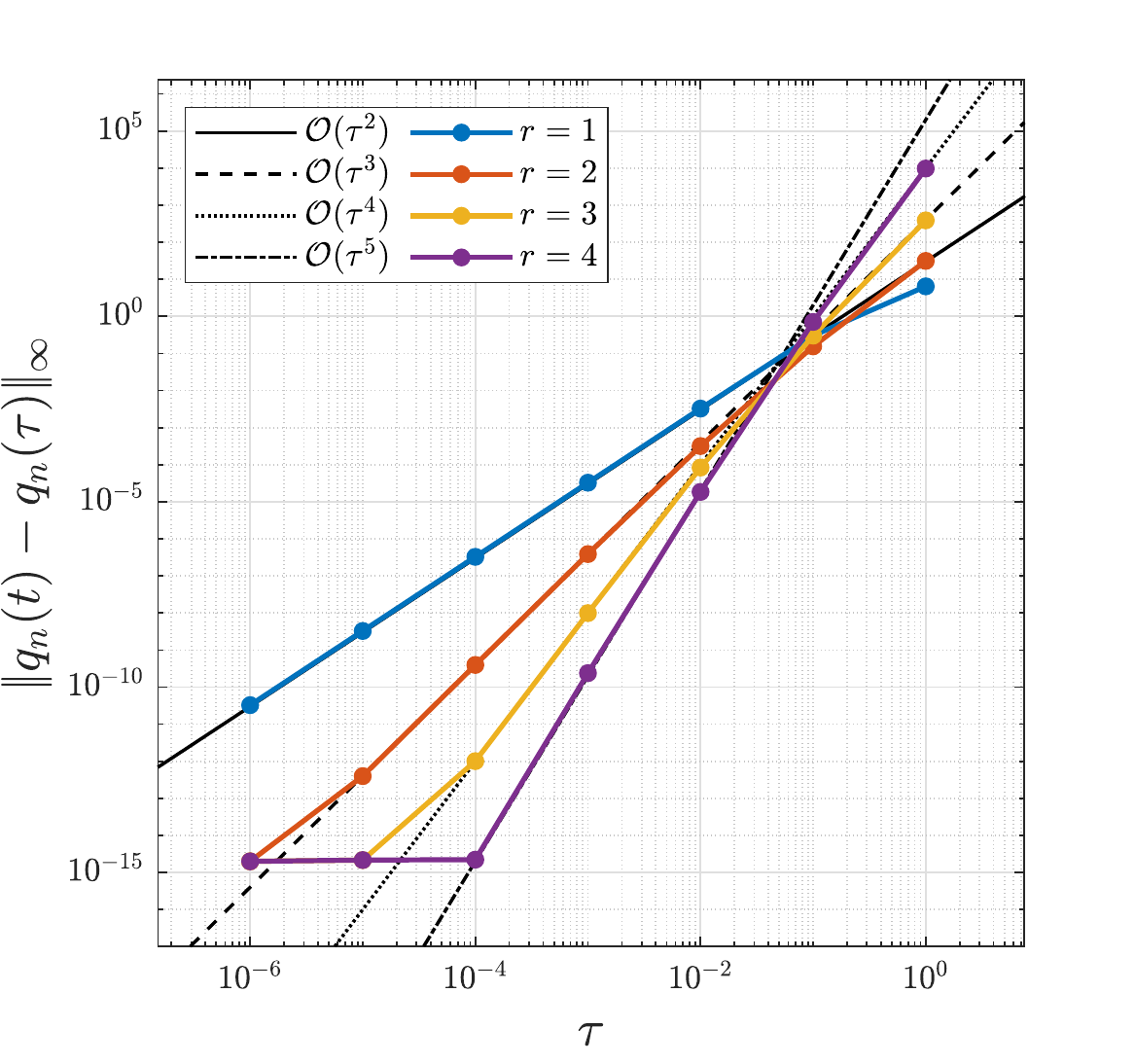}
  		\caption{}
  		\label{advec2_splitstep_b}
	\end{subfigure}	
	\caption{(a) The error (green) between the SD-UTM solution \eqref{soln_advec2_backward} and the exact solution and the error (maroon) between the exact solution and the small-time approximated SD-UTM solution \eqref{soln_advec2_splitstep} with $r = 2$ and $t = \tau = 10^{-5}$ as $h \rightarrow 0$. (b) The error between the SD-UTM solution \eqref{soln_advec2_backward} and the small-time approximated SD-UTM solution \eqref{soln_advec2_splitstep} for varying $r$ as $\tau \rightarrow 0$ with $h = 0.01$.}
	\label{advec2_splitstep}
	\end{figure}

	As an example, consider the IBVP
	\begin{equation}
	\begin{dcases}
		q_t = - q_{x},& x > 0,\, t > t_0, \\
		q\left(x,t_0\right) = \phi(x) = \frac{e^{-2x} \left( \sin \left(4 \pi x \right) + 1\right)}{2},& x > 0, \\
		q(0,t) = u(t) = \frac{1}{2} + (1 - 2\pi)t e^{-t},& t > t_0,
	\end{dcases}
	\label{advec2_numerical1_HL}
	\end{equation}	
	with $t_0 = 0$. Figure \ref{advec2_splitstep_a} depicts the errors, relative to the exact solution, as $h \rightarrow 0$ for the semi-discrete UTM solution \eqref{soln_advec2_backward} and \textit{small-time approximated} semi-discrete UTM solution \eqref{soln_advec2_splitstep} up to terms of order $2$. Despite only including terms up to second order, the errors in the plot are indistinguishable, implying we need not include higher-order terms to obtain a suitable small-time approximate solution \eqref{soln_advec2_splitstep}. Figure \ref{advec2_splitstep_b} depicts the errors for the small-time solution \eqref{soln_advec2_splitstep} as $\tau \rightarrow 0$ for a fixed $h$, relative to \eqref{soln_advec2_backward}. This plot shows that relatively large values of $\tau$ lead to an accurate approximation to \eqref{soln_advec2_backward}.


\section{Concluding Remarks}

	The SD-UTM is applied algorithmically using the following steps:
\begin{enumerate}
	\item Rewrite into divergence form to obtain the local relation and the dispersion relation $W(k)$, \label{sd_utm_1}
	\item Sum over spatial indices and integrate over the temporal domain to obtain the global relation, \label{sd_utm_2}
	\item Invert to obtain a representation of the solution depending on unspecified boundary data, \label{sd_utm_3}
	\item Determine the symmetries $\nu_j(k)$ of $W(k)$, \label{sd_utm_4}
	\item Determine where the global relations with $k \rightarrow \nu_j(k)$ are valid, \label{sd_utm_5}
	\item If necessary, deform integral paths of the boundary terms appropriately, \label{sd_utm_6}
	\item \textit{If necessary, determine additional boundary conditions from the PDE}, \label{sd_utm_7}
	\item \textit{Appropriately discretize boundary conditions}, \label{sd_utm_8}
	\item Solve for unknowns using global relations with $k \rightarrow \nu_j(k)$ and time transforms of discretized boundary conditions, and \label{sd_utm_9}
	\item Check integral terms involving $\hat{q}(\nu_j,T)$ vanish, resulting in the solution representation depending only on known quantities. \label{sd_utm_10}
\end{enumerate}
	With minor differences in the calculations, the procedure for the semi-discrete UTM is almost identical to that from the continuous UTM, with Steps \eqref{sd_utm_7} and \eqref{sd_utm_8} added. The SD-UTM operates similarly for finite-interval IBVPs, presented in a forthcoming paper \cite{SDUTM_FI}. The steps themselves become more tedious for higher-order problems, but like the continuous UTM, the SD-UTM reduces the burden of solving a semi-discrete IBVP to solving for the roots of polynomials and dealing with a set of algebraic equations. Third-order problems, like the linear Korteweg-de Vries equations $q_t = \pm q_{xxx}$, will be presented in a forthcoming paper \cite{SDUTM_kdv}.
	
	With the SD-UTM, we develop the notion of ``natural'' discretizations, which reduce the variety of stencils down to those that are compatible with the IBVP. The natural discretization is (i) of the same order as the spatial order of the PDE, (ii) not purely one sided (except for first-order problems), and (iii) the one that optimally aligns with the available boundary conditions. Once the PDE is discretized, it follows that the available discretizations for any derivative boundary conditions are dictated by the global relation and its validity under the symmetries $\nu_j(k)$, as we saw with the Neumann IBVP \eqref{heat_prob_N} and \eqref{LS_prob_N}.


\section{Acknowledgments}
	
	This work was supported by the Graduate Opportunities \& Minority Achievement Program Fellowship from the University of Washington and the Ford Foundation Predoctoral Fellowship (JC). Any opinions, findings, and conclusions or recommendations expressed in this material are those of the authors and do not necessarily reflect the views of the funding sources.


\vspace{10pt}
\nocite{*}
{\small
\bibliographystyle{abbrv}
\bibliography{references}
}

\end{document}